\documentclass[12pt,dvips]{amsart}

\title{Complete moduli spaces of branchvarieties}

\newif\iffancyformatting
\fancyformattingtrue
\usepackage{amsfonts, amssymb, latexsym, epsfig,epic}

\iffancyformatting
\usepackage{euler}
\fi

\author{Valery Alexeev}
\address{Department of Mathematics, University of Georgia, 
  Athens GA 30602, USA}
\email{valery@math.uga.edu}
\author{Allen Knutson}
\address{Department of Mathematics, UCSD, La Jolla CA 92093, USA}
\email{allenk@math.ucsd.edu}
\date{July 31, 2006}

\iffancyformatting

\fi

\newcommand{\ol}{\overline}
\newcommand{\mb}{\mathbb}
\newcommand{\mc}{\mathcal}
\newcommand{\wt}{\widetilde}

\newcommand{\inv}{^{-1}}
\newcommand{\sat}{^{\rm sat}}
\newcommand{\red}{_{\rm red}}
\newcommand\cO{{\mathcal O}}
\newcommand\bP{{\mathbb P}}
\newcommand\bA{{\mathbb A}}

\DeclareMathOperator{\lcm}{\operatorname{lcm}}

\newcommand{\Chow}{\operatorname{Chow}}
\newcommand{\Forest}{\operatorname{Forest}}
\newcommand{\PGL}{\operatorname{PGL}}
\newcommand{\Pic}{\operatorname{Pic}}
\newcommand{\Aut}{\operatorname{Aut}}
\newcommand{\Ext}{\operatorname{Ext}}
\newcommand{\Hom}{\operatorname{Hom}}
\newcommand{\Hilb}{\operatorname{Hilb}}
\newcommand{\B}{\operatorname{Branch}}
\newcommand{\Bhb}{\operatorname{Branch}_{h,Y}^b}
\newcommand{\cB}{\mathit{Branch}}

\newcommand{\Spec}{\operatorname{Spec}}
\newcommand{\Proj}{\operatorname{Proj}}

\newcommand{\codim}{\operatorname{codim}}
\newcommand{\im}{\operatorname{im}}
\newcommand{\chr}{\operatorname{char}}
\newcommand{\diag}{\operatorname{diag}}

\DeclareMathSymbol{\twoheadrightarrow}  {\mathrel}{AMSa}{"10}
\newcommand{\onto}{\twoheadrightarrow}

\newtheorem{theorem}{Theorem}[section]
\newtheorem{lemma}[theorem]{Lemma}
\newtheorem{proposition}[theorem]{Proposition}
\newtheorem{corollary}[theorem]{Corollary}
\theoremstyle{definition}
\newtheorem{definition}[theorem]{Definition}
\newtheorem{assumption}[theorem]{Assumption}
\newtheorem{example}[theorem]{Example}
\newtheorem{question}[theorem]{Question}
\newtheorem{remark}[theorem]{Remark}
\newtheorem*{acknowledgements}{Acknowledgements}
\newcommand\op{{\rm opp}}

\newcommand\PP{{\mathbb P}}
\renewcommand\AA{{\mathbb A}}
\newcommand\tensor{\otimes}
\newcommand\naturals{{\mathbb N}}
\newcommand\defn[1]{{\bf #1}} % to allow for \emph{#1}
\newcommand\junk[1]{}
\newcommand\into{\hookrightarrow}

\begin{document}

\pagestyle{plain}
\begin{abstract}
  The space of subvarieties of $\PP^n$ with a fixed Hilbert polynomial
  is not complete. Grothendieck defined a completion by relaxing
  ``variety'' to ``scheme'', giving the complete \emph{Hilbert scheme}
  of subschemes of $\PP^n$ with fixed Hilbert polynomial.
  
  We instead relax ``sub'' to ``branch'', where a \defn{branchvariety of} 
  $\PP^n$ is defined to be a \emph{reduced} (though possibly reducible)
  scheme \emph{with a finite morphism to} $\PP^n$. Our main theorems are that
  the moduli stack of branchvarieties of $\PP^n$ with fixed Hilbert
  polynomial and total degrees of $i$-dimensional components is a 
  proper (complete and separated) Artin stack with finite stabilizer,
  and has a coarse moduli space which is a proper algebraic space.
  
  Families of branchvarieties have many more locally constant invariants than
  families of subschemes; for example, the number of connected components
  is a new invariant. In characteristic $0$, one can extend this count
  to associate a $\mb Z$-labeled rooted forest to any branchvariety.
\end{abstract}
\maketitle

\iffancyformatting{%\small
\baselineskip 10pt
\tableofcontents
}\else{%
\tableofcontents
}\fi

\setcounter{section}{-1}

\section{Introduction}

Consider a family of reduced plane conics $\{ x_1^2 = t\cdot x_0x_2 \}$ 
over $\mb A^1_t\setminus 0$ with coordinate~$t$.   By using the same
equation, this family can be completed to a flat family $Z\subset \mb P^2
\times \mb A^1_t$ projective over $\bA^1_t$; the central fiber is
the double line $Z_0 = \{x_1^2=0\}$, a nonreduced subscheme of $\mb P^2$.
On the other hand, after the finite ramified base change $t=s^2$ and
simple substitution $x'_1 = x_1/s$, this family can be completed to a
flat family $X=\{ (x_1')^2 = x_0x_2 \}$ with a \emph{finite morphism}
$f:X\to \mb P^2 \times \mb A^1_s$, rather than an inclusion. The
central fiber $X_0$ is a \emph{reduced} $\mb P^1$, and $f_0:X_0\to \bP^2$ 
is a double cover of the line $\{x_1=0\} = (Z_0)\red$.

How general is this phenomenon? We pose two forms of this question:
\begin{enumerate}
\item Does every one-parameter flat family %over $\mb A^1_t\setminus 0$
  of reduced subschemes of $\PP^n$ have a well-defined flat completion 
  whose special fiber is a reduced scheme mapping finitely to $\PP^n$?
  (Possibly after base change, like $s$ for $t$ in the above example.)
\item Is there a universal substitute for the Hilbert scheme, some
  other moduli space for reduced schemes $X$ carrying finite
  morphisms $X\to \mb P^n$? (Such a moduli space would necessarily be
  coarse since $\Aut(f)$ may be nontrivial; indeed it is $\mb Z_2$ for
  the double cover $f_0$ above.)
\end{enumerate}

In this paper we answer both of these positively,
generalizing two known situations:

\begin{enumerate}
\item In \cite{Alexeev02,AlexeevBrion05} there were constructed the moduli
  spaces of stable toric and spherical varieties \emph{over $Y$}. This can be
  interpreted as answering both questions above in the \emph{multigraded
    multiplicity-free case}.  The one-parameter limits were constructed by
  making a base change $t=s^n$ and normalizing, which in the multiplicity-free
  case amounts to saturating a semigroup. These moduli spaces are projective,
  and they should be considered to be the substitute for the toric Hilbert
  scheme of \cite{PeevaStillman}, the multigraded multiplicity-free Hilbert
  scheme.

\item For every closed subscheme $X$ of a reduced Noetherian scheme $Y$,
  there is a flat degeneration of $Y$ to the {\em normal cone} $C_X Y$,
  which is usually nonreduced. In \cite{Knutson05} there is presented
  an alternative degeneration of $Y$ to a ``balanced'' normal cone
  $\ol C_X Y$ which {\em is} reduced, and comes with a finite morphism 
  $\ol C_XY\to C_XY$. 
The schemes $X,Y$ are otherwise general, and there is no multiplicity-free
restriction or group action.
\end{enumerate}

The need for the base change, and the loosening of ``inclusion of a
subscheme'' to ``finite morphism'', are already both illustrated by
the simple example $f:\AA^1_u \to \AA^1_t$, $t = u^2$. The general
fiber of this flat family is two (reduced) points, which cannot be
disentangled because the source is irreducible. Pulling back this family
along the base change $s^2 = t$, we get another flat family
$f':\AA^1_{s=u} \cup_0 \AA^1_{s=-u} \to \AA^1_s$ with the same central fiber 
(a double point), where the total space of the family is given by $s^2=u^2$.
Then we normalize to pull the two components apart; now the central
fiber too is two reduced points.

A similar normalization step will introduce finite morphisms in general, 
taking us out of the class of subvarieties (i.e. {\em injective} morphisms of
reduced schemes). We now spell out what class will replace them.

\begin{definition}
  Let $Y\subset \bP^n$ be a projective scheme over a field $k$. 
  A \defn{branchvariety of $\mathbf{Y}$} is a variety $X$ (by which we mean a
  scheme of finite type over $k$ such that $\ol X=X\otimes_k \bar k$
  is reduced, so \emph{a forteriori}, $X$ is reduced) equipped with a
  finite (hence proper) morphism $f:X\to Y$.
  
  The \defn{Hilbert polynomial} of a branchvariety is $h(d) := \chi(X,L^d)$,
  where $L= f^*\mc O_{\bP^n}(1)$.  
  One has $\dim H^0(X,L^d)=h(d)$ for $d\gg0$.
  Writing the leading term of $h(d)$ as $(c/D!)\, d^D$, the coefficient
  $c$ is the \defn{degree} $\deg X$ of the branchvariety.
%  The \defn{automorphism group} of $f:X\to Y$ is the automorphisms $g:X\to X$
%  such that $f\circ g = f$.
\end{definition}

We could even assume without much loss of generality that $X$ is
connected, since by Lemma \ref{lem:connectedcomponents} to come (and
in contrast to families of subschemes) the number of connected
components is locally constant in families of branchvarieties.
Then $X$ is indeed what is commonly called a ``variety'', albeit
not irreducible, much as a ``stable curve'' is not irreducible in general.
We also note that one can define in the same way a branchvariety $X$
of an arbitrary, not necessarily projective scheme $Y$, though giving
up the notion of $X$'s Hilbert polynomial and degree.

Let us fix a Noetherian base scheme $C$ (for example the spectrum
of $\mb Z$ or of a field of arbitrary characteristic), and a projective scheme
$Y\subset \bP^n_C$ flat over $C$. 
All our schemes will be assumed to be locally Noetherian (but see
Remark~\ref{non-Noetherian}).

\begin{definition}
  The \defn{functor of subschemes} or \defn{Hilbert functor}
  $$\Hilb_{h,Y}: \text{($C$-schemes)} \to \text{(Sets)}^\op$$
  is defined by associating to each 
  scheme $S$ the set $\Hilb_{h,Y}(S)$
  of proper subschemes $Z\subset Y_S= Y\times_C S$ which are flat over $S$.
\end{definition}

Because of the importance of flatness in this and the next definition,
in this paper we will only be interested in flat families, and will
{\em always} use the term ``family'' to mean ``flat family''.

One fixes the Hilbert polynomial in the above definition in order that
the functor be representable by a scheme of finite type. (Otherwise
one gets a disjoint union of schemes of finite type, one for each possible
Hilbert polynomial.) In the corresponding definition for branchvarieties,
fixing the Hilbert polynomial is not enough to obtain a finite-type family, 
as Example~\ref{ex:unbounded} will show. One easy workaround is to
look at only equidimensional branchvarieties; instead, to avoid loss of
generality, we measure some additional parameters beyond $h$. 
For a branchvariety $X$, let
$$ X^{\dim i} := \bigcup\ \left\{ 
  \text{the $i$-dimensional irreducible components of $X$}\right\}, \qquad
  b_i(X) := \deg X^{\dim i}$$ 
where $X^{\dim i}$ is considered as a
branchvariety in a natural way and hence has a degree.  
By Lemma~\ref{lem:degree-sequence} the
\defn{degree sequence} $b=(b_0,\dotsc, b_{\dim X})$ of nonnegative integers
%where $d=\dim X=\deg h$, 
is locally constant in families of branchvarieties.

\begin{definition}
  The \defn{functor of branchvarieties} 
  $$\B^b_{h,Y}: \text{($C$-schemes)} \to \text{(Sets)}^\op$$
  is defined by
  associating to each scheme $S$, the set $\B_{h,Y}(S)$ of proper families
  $f:X\to Y_S = Y\times_C S$ such that $X\to S$ is flat and every fiber
  $f_s:X_s\to Y_s=Y\times_C k(s)$ is a branchvariety of $Y_s$ with Hilbert
  polynomial $h$ and degree sequence $b$, up to isomorphism.

%  such that for $d\gg 0$ the sheaves
%  $(p_2\circ f)_*(L^d)$ are locally free sheaves on $S$ of rank $h(d)$, 
%  where $L = f^*\mc O_{\mb P^n_S}(1)$, 

  We avoid the set-theoretic difficulties in this definition (the
  ``set'' of branchvarieties may be too big to actually be a set) 
  in a standard way, by fixing a universe, or by demanding that they
  are subschemes of some fixed $\PP^\infty$. We prove in Theorem
  \ref{thm:boundedness} that $\infty$ can actually be lowered to a
  finite bound $d_0(h,b)$.
  
  We define the \defn{stack in groupoids} $\cB^b_{h,Y}$ by associating to each
  scheme $S$ a category $\cB^b_{h,Y}(S)$ whose objects are the families
  $f:X\to Y_S$ as above, and morphisms are isomorphisms
  intertwining the structure morphisms $f$.
\end{definition}

\begin{theorem}\label{thm:main}
  The stack in groupoids $\cB^b_{h,Y}$ is an algebraic Artin stack with a
  finite stabilizer. It has a coarse moduli space $\Bhb$ which is a proper
  algebraic space.
\end{theorem}

In particular, $\Bhb$ has finitely many connected components. 
It is not connected for most $(h,b,Y)$, 
since different connected components may often be distinguished
by their associated (labeled rooted) ``forests'',
an invariant we develop in Section \ref{sec:rooted-forests}.
(This is in contrast with the Hilbert scheme, which is connected for
each $h$.)
The forest can be defined when the characteristic of our field is $0$ or
larger than any degree $b_i(X)$, and is a refinement of both the
degree sequence and the Hilbert polynomial.

The characterization of the set of connected components of $\Bhb$
remains open, though we have one tiny result in this direction in
Theorem \ref{thm:Kollar}.

\medskip

The proof of Theorem~\ref{thm:main} uses the general procedure which was
developed for constructing compactifications of moduli spaces of surfaces of
general type (moduli of stable surfaces and pairs), see e.g.
\cite{Kollar_Projectivity,Alexeev_Mgn}, and it goes as follows.

In Section \ref{sec:one-param-families} we establish the most important
property of our moduli functor, properness: every one-parameter family
of branchvarieties has at most one limit, and the limit always exists after a
finite ramified base change.

In Section \ref{sec:boundedness} we prove that the family of branchvarieties
with fixed numerical invariants is bounded, i.e. there exists a universal
constant $N$ such that $L^N$ is very ample for each such branchvariety and
embeds it into a fixed projective space $\bP^D$.  Once this is established,
the branchvarieties of $Y$ can be parametrized by a locally closed subscheme
$V$ of the Hilbert scheme of $Y\times \bP^D$, up to a choice of embedding
$X\into\bP^D$. The moduli space can then be constructed by taking the quotient
$V/\PGL_{D+1}$. We do this in Section \ref{sec:constr-moduli-space}.

For the quotient, we do not use Geometric Invariant Theory, which would have
involved a delicate analysis of stability. Instead, we use a well-known
observation that a quotient by a proper group action always exists as an
algebraic space.  Properness of the group action follows from 
separatedness of the moduli functor and finiteness of the automorphism
groups.  Since the moduli functor of branchvarieties is also proper, the
moduli space is a proper algebraic space.

In the case of stable surfaces of general type, there is one additional step
one can add to 
this procedure: by using \cite{Kollar_Projectivity} one proves that the thus
obtained moduli space is projective, and in particular a scheme, because some
naturally defined invertible sheaves on it are ample. In our case this part is
missing.  As we show in Section~\ref{sec:line-bundles-moduli}, the basic
invertible sheaves that come with $\B$ are, curiously, not ample. However, it
is possible that some of their linear combinations are ample.

\begin{remark}\label{rem:Branch-vs-Chow}
  We would like to note that Branch does not suffer from
  the limitations afflicting the Chow variety (and the similar Barlet
  space in complex-analytic geometry) which parametrizes cycles of a
  fixed degree in $Y$. In the Chow theory, there is no notion of an
  infinitesimal family of algebraic cycles, say over
  $\Spec \mb C[\epsilon]/(\epsilon^2)$; the best one can do is to
  consider families of algebraic cycles over a reduced and seminormal
  base.  In this sense, the Chow variety should be considered to
  be a ``parameter space'' rather than a ``moduli space''.
  See \cite[I.3-4]{Kollar_RationalCurves} for a comprehensive rigorous
  treatment. 
\end{remark}

\begin{remark}\label{non-Noetherian}
  The Hilbert scheme has been constructed in the more general situation of
  families over non-Noetherian schemes.  In this case, the Hilbert functor is
  defined by requiring that the sheaves $(p_2)_*\mc O_Z(d)$ are locally free
  of rank $h(d)$, in a neighborhood of every point $s\in S$ and for $d\ge
  d_0(s)$.  This is equivalent to flatness in the Noetherian case since a
  finitely generated module over a Noetherian ring is flat iff it is locally
  free.
  
 Similarly, we can define a family of branchvarieties over a general scheme $S$
  by requiring that the sheaves $\pi_* L^d$ are locally free of rank $h(d)$.
  The moduli space of branchvarieties is constructed in Theorem \ref{thm:main}
  by taking a quotient of a Hilbert scheme, which is already quasiprojective
  over $C$. Then, in order to use the result from \cite{KeelMori}
  concerning this quotient, we only need to assume that $C$ is locally
  Noetherian.
\end{remark}

\begin{remark}
  Since the moment this paper was widely circulated, two extensions have
  already appeared: M. Lieblich \cite{Lieblich} and J. Starr \cite{Starr}
  extended some of our constructions to the case of an Artin stack $Y$ as the
  target.  In addition, \cite{Alexeev_Limits} contains an application of
  branchvarieties to the moduli of weighted stable pairs.
\end{remark}

\begin{acknowledgements}
  We would like to thank H\'el\`ene Esnault for an argument with Frobenius
  used in the proof of Theorem~\ref{thm:boundedness}, and J\'anos Koll\'ar for
  suggesting Theorem~\ref{thm:Kollar}; of course any errors remaining are
  ours. Also, we thank Patricia Hersh for discussions about rooted forests,
  David Speyer for correcting our calculation of stable cubics,
  and Diane Maclagan for asking about multigraded $b$-sheaves.
  
  Michael Thaddeus informed us that he proposed the functor of this
  paper in the case of curves in a 1995 talk at Harvard, with the
  family of plane cubics as evidence (treated here in Section 
  \ref{sec:complete-cubics}). We also note that Morten H\o nsen
  constructed a proper moduli space of curves $X$ with a finite morphism
  $f:X\to \mb P^n$ \cite{Honsen}; in place of our assumption that $X$
  is reduced ($S1$ and $R0$), he only requires $X$ to be $S1$,
  but also that $f$ be generically an embedding.
\end{acknowledgements}

%%%%%%%%%%%%%%%%%%%%%%%%%%%%%%%%%%%%%%%%
\section{Boundedness}
\label{sec:boundedness}

We will now fix a certain class of branchvarieties, and show that this
class is \defn{bounded}. By this we will mean that for every
branchvariety $f:X\to \bP^n$ in our class,
over an algebraically closed field $k=\bar k$, 
some \emph{fixed} power $f^* \mc O_{\bP^n}(1)^{d_0}$ is very ample.
Eventually this will give us a dimension $D$ to use in the
construction of the moduli space described in the Introduction.

\begin{example}\label{ex:unbounded}
  Let $n$ be a nonnegative integer and $X$ be a union of two $\mb P^1$s
  joined at $n$ simple nodes, %$P_1, \dotsc, P_n$ 
  plus a disjoint union of $n$ points. % $Q_1,\dotsc, Q_n$. 
  For any such $(X,f:X\to \mb P^1$ of degree $2)$, 
  the Hilbert polynomial is the same, $h(d)=2(d+1)$. 
  But as Lemma \ref{lem:degree-sequence} will show, if the degree sequences
  $b_0=n$, $b_1=2$ of these curves were to all appear in a single family,
  the base scheme would need infinitely many connected components and
  hence not be of finite type.
  
  This shows that branchvarieties with a fixed Hilbert polynomial do not form
  a bounded class unless one adds some further restrictions. The easiest would
  be to require that $X$ be equidimensional.  We describe a more general
  solution below, which we refine further in Section \ref{sec:rooted-forests}
  in the case $\chr k=0$ or large enough.
\end{example}

Let us introduce some notation. Let $Z\subset \mb P^n$ be the image of $X$,
with reduced scheme structure. Choose a sequence of general linear forms $l_1,
\dots, l_{\dim X}$ on $\mb P^n$, and define $Z_i, X_i$ inductively as the
Cartier divisors $l_i=0$ in $Z_{i-1}, X_{i-1}$ respectively, with $Z_0=Z$,
$X_0=X$. We know that the $Z_i$ are reduced irrespective of $\chr k$
(\cite[II.8.18]{Hartshorne}, \cite{Flenner}) but $X_i$ need not be reduced if
$f:X\to Z$ is not separated: consider for example the geometric Frobenius map
$f:\mb P^1\to \mb P^1$,  $x\mapsto x^p$, in which $X_1$ is seen to be 
a point of multiplicity $p$. 

\begin{remark}\label{rem:big-field}
  The precise generality condition on the linear forms is that each
  $l_i$ does not vanish identically on any of a certain finite set of proper
  subvarieties of $\mb P^n$, namely, the associated components of
  $Z_{i-1}$ and the subvarieties appearing in Bertini's theorem. For
  the proofs of the statements of this Section we are free to make
  finite base extensions, and may thereby assume that such general linear
  forms do indeed exist.
\end{remark}

\begin{lemma}\label{lem:degree-sequence}
  The integers $b_i=\deg X^{\dim i}$ are locally constant in
  families of branchvarieties.
\end{lemma}

\begin{proof}
  By making a base change $\Spec A \to S$, it is sufficient to consider the
  case of a regular one-dimensional base, for example $A$ a DVR. Any
  Cartier divisor $l_1=0$ on the central fiber $X_0$ can be extended
  to a Cartier divisor on the generic fiber $X_{\eta}$. Hence, the Cartier
  divisors $X_1$ form a flat family over $\Spec A$ with the Hilbert
  polynomial $p_1(d)=p(d)-p(d-1)$. By induction we see that all $X_i$ can be
  put in flat families.

  Now, $\deg X^{\dim i}= \deg (X_i)^{\dim 0}$. The latter space is a union of
  connected components of $X_i$. Hence, it too is flat over $\Spec A$, and
  its length (or cardinality, if it is reduced and $k=\bar k$) is constant.  

\end{proof}

Recall the following definitions:

\begin{definition}[Kleiman]
  Let $b=(b_0,b_1, \dotsc, b_n)$ be a sequence of integers. A coherent sheaf
  $F$ on $\bP^n$ is a \defn{$\mathbf b$-sheaf} if for generic hyperplanes $l_1,
  \dotsc, l_n$ and the inductively defined sheaves 
  $F_0=F$, $F_i= F_{i-1}/l_i F_{i-1}$, one has $h^0( F_i(-1)) \le b_i$.
\end{definition}

\begin{definition}
  A coherent sheaf $F$ on $\bP^n$ is said to be Castelnuovo-Mumford
  \defn{$\mathbf m$-regular} if $H^i( F(m-i))=0$ for all $i>0$.
\end{definition}

The following is possibly the strongest known result implying boundedness of
various classes of coherent sheaves.

\begin{theorem}[Kleiman, \cite{SGA6}, Thm. XIII.1.11]
  \label{thm:Kleiman}
  For fixed $b$ and $h(d)$ there exists an integer $m$ such that every
  $b$-sheaf $F$ with Hilbert polynomial $\chi(F(d))=h(d)$ is $m$-regular.
\end{theorem}

Here is our application of this result.

\begin{theorem}\label{thm:boundedness}
  Fix a Hilbert polynomial $h(d)$ and nonnegative integers $b_0, b_1,
  \dotsc, b_{\deg h}$. Then there exists a positive integer $d_0$ such that
  the following holds:
  
  For any branchvariety $f:X\to \bP^n$ with Hilbert polynomial $h(d)$ and
  $\deg X^{\dim i}=b_i\ \forall i$, the sheaf $L^d$ is very ample for
  $d\ge d_0$ and the algebra
  $$ R(X,L)^{(d_0)} = k \oplus \bigoplus_{d\ge 1} H^0(X,L^{dd_0})$$
  is generated in degree $1$. 
\end{theorem}

(Although we won't need it, for a slightly larger $d_0$ one can also
ensure that the relations are generated in degree $2$.)

\begin{proof}
  By \cite[Thm. 3]{Mumford_QE}, the statement would follow if we could prove
  that $H^i(X,L^d)=0$ for all $d\ge d_1$ and $i>0$ (we note the importance of
  the fact that $L$ is free, which we have): then $L^{(\dim X+1)d_1}$ is very
  ample. 

  We apply Theorem~\ref{thm:Kleiman} to the sheaf $F=f_*\mc O_X$. Then $F_i=
  f_*(\mc O_{X_i})$ and $H^0( F_i(-1) ) = H^0( X_i, L^{\inv})$.
  If $\chr k=0$, the generic sections $X_i$ are reduced by Bertini's theorem,
  so $h^0( X_i, L^{\inv}) = \deg X_i^{\dim 0} = b_i$, and we are done. 
  
  If $\chr k = p>0$, one has to be a little more careful. Decompose $X\to Z$
  into a purely inseparable morphism $X\to Y$ followed by a separable morphism
  $Y\to Z$. Then all $X_i\to Z_i$ decompose into purely inseparable $X_i\to
  Y_i$ followed by separable $Y_i\to Z_i$, and for generic hyperplanes $l_i$
  the schemes $Y_i$ are reduced.  
  Since $X_i\to Y_i$ is purely inseparable,
  it can be dominated $F^k(Y_i)\to X_i \to Y_i$
  by a power $F^k(Y_i)$ of the absolute Frobenius $F(Y_i)$.
  Hence,
  $$ H^0( \mc O_{X_{i+1}} ) \subset H^0( \mc O_{ p^k Y_{i+1}} ), $$
  where the latter is the Cartier divisor on $Y_i$ with the equation
  $l_{i+1}^{p^k}$.
  
  Now, let $D$ be a reduced ample Cartier divisor on a projective scheme $Y$
  and assume that every connected component of $D$ has dimension $\ge 1$. Then
  the basic exact sequence
  $$ 0\to \mc O_D(-(s-1)D)\to \mc O_{sD} \to \mc O_{(s-1)D} \to 0 $$
  implies that 
  $$ h^0 (\mc O_{sD}(-D)) \le h^0 (\mc O_{(s-1) D}(-D)) \le \dotsb
  \le h^0 (\mc O_{D}(-D))=0 $$
  
  Applying this to $D=Y_{i+1}$ on $Y=Y_i$, we see that again only the
  zero-dimensional connected components contribute, and so 
  $h^0(X_i, L^{-1}) = \deg X_i^{\dim 0} = b_i$ as before.
\end{proof}

\begin{theorem}\label{thm:family=algebra}
  Fix a Hilbert polynomial $h$ and degree sequence $b$, and let $d_0$
  be as in Theorem \ref{thm:boundedness}.
  Then a family of branchvarieties over $Y_{\Spec A}$ 
  with Hilbert polynomial $h(d)$ and with each $\deg X^{\dim i}=b_i$ 
  is equivalent to a graded ring $R$ together with a homomorphism
  $\phi:A[x_0, \dotsc, x_N]\to R$, $N= \binom {n+d_0} {n}$ such that
  \begin{enumerate}
  \item $R_0=A$, 
  \item the $R_d$ are finite locally free $A$-modules of rank $h(dd_0)$,
  \item $R$ is finite over the image of $\phi$,
  \item $\ker \phi$ contains the ideal of the $d_0$-tuple Veronese image of
    $Y$, and
  \item for every homomorphism $A\to \bar k$ to an algebraically closed
    field, $R\otimes_{A} \bar k$ is reduced.
  \end{enumerate}
  Moreover, any such algebra $R$ is generated over $A$ in degree 1.
\end{theorem}
\begin{proof}
  We first note that under the isomorphism between $Y$ and its Veronese
  embedding $v_{d_0}(Y)$, the families of branchvarieties $f_1:X\to Y_S$ and
  $f_2:X\to v_{d_0}(Y)_S$ are in natural bijective correspondence. The
  (same) invertible sheaf $L$ is uniquely determined by either $f_1$ or $f_2$.
  Hence, without loss of generality we can replace $Y$ by $Y'=v_{d_0}(Y)$.

  Given a family $f:X\to Y'_S$, we set 
  $$
  R= R(X/S,L')^{(d_0)} = A \oplus \bigoplus_{d\ge 1} H^0(X,(L')^d),
  \quad
  \text{where } L'=f_2^*\mc O_{\bP^N}(1) = L^{d_0}.
  $$
  Since the higher cohomology groups of $(L')^d$ vanish, by the
  Cohomology and Base Change Theorem $H^0(X,(L')^k)$ are locally free
  modules of rank $h(dd_0)$.

  In the opposite direction, we set $X=\Proj R$, and the condition (3) gives a
  morphism $f:X\to \bP^N_A$, which factors through $Y'$ by the condition
  (4). Clearly, the associations $(R,\phi)\leftrightarrow (X,f)$ 
  are inverses of each other, and the condition (5) is equivalent to
  the condition that the geometric fibers of $X$ are reduced.

  The last statement is a direct consequence of Theorem \ref{thm:boundedness}.
\end{proof}

%%%%%%%%%%%%%%%%%%%%%%%%%%%%%%%%%%%%%%%%
\section{One-parameter families}
\label{sec:one-param-families}

In this section, $A$ is a DVR with maximal ideal $(t)$, residue field
$k=A/tA$, and fraction field $K=A[1/t]$. 
The ring $R=\oplus_{d\ge0} R_d$ is a finitely generated $A$-algebra, and
each $R_d$ is a locally free $A$-module of rank $h(d) < \infty$. 
By Theorem \ref{thm:family=algebra}, after replacing $Y\subset \bP^n$ with its
Veronese embedding $v_{d_0}(Y)\subset \bP^N$, such an algebra $R$ is
equivalent to a family of branchvarieties over $\Spec A$.

Let $\wt R$ denote the integral closure of $R$ in $R_K := R\tensor_A K$, 
i.e. in the \defn{general fiber} of the family $\Spec R \to \Spec A$.
Since by definition $\wt R \subset \oplus_{d\ge0}\ R_d\otimes_A K$, it
is also a graded ring. A ring finitely generated over a DVR has a
finite normalization, and therefore $\wt R$ is a finitely generated
$R$-module. It follows that each $\wt R_d$ is again a finitely
generated $A$-module, so it is free of the same rank $h(d)$ and
for each $d$, there exists a $k_d$ such that 
$R_d \subset \wt R_d \subset  t^{-k_d} R_d$.

By an \defn{$\mathbf n$-ramified base change} we will mean another DVR 
$A'\supset A$ with uniformizing
parameter $t'$, so that $t= c'(t')^n$ for some unit $c'\in A'$.
We will use $R'$ for $R\otimes_A A'$, $K'$ for the fraction field,
and $\wt R'$ for the integral closure $\widetilde{(R')}$ of $R'$ in~$K'$.

\begin{lemma}\label{lem:nilpotents}
  \begin{enumerate}
  \item $R/(t)$ is reduced $\implies$ $\wt R= R$.
  \item $R/(t)$ is not reduced $\implies$ after some ramified base change $\wt
    R' \ne R'$.
  \end{enumerate}
\end{lemma}
\begin{proof}
  (1) Suppose $\wt R\ne R$. Then there exists $x\in R\setminus tR$ such that
  $x/t \in 
  \wt R \setminus R$, and it satisfies some monic equation
  \begin{eqnarray*}
    \left(\frac x t\right)^n + r_{n-1} \left(\frac x t\right)^{n-1} 
    + \dotsb + r_0 = 0  
    &\implies&
     x^n + tr_{n-1} x^{n-1}  + \dotsb + t^n r_0 = 0  
     \\
    &\implies& x^n \in (t) \\ 
    &\implies& \bar x\in R/(t)
    \text{ is a nonzero nilpotent.}
  \end{eqnarray*}
 
  (2) Suppose some $\bar x\in R/(t)$ is a nonzero nilpotent, i.e.
  there exists some $x\in R \setminus tR$ such that $x^n\in tR$. Then
  after an $n$-ramified base change one has $x^n\in (t')^n R'$. 
  So $(x/t')^n \in R'$. Hence, $x/t'$ is integral over $R'$, and so 
  $x/t' \in \wt R'$. It remains to show $x/t' \notin R'$.

  Since $R,R'$ are torsion-free $A,A'$-modules over DVRs, for any $x\in R$
  we have the equivalences
  $$ x/t \in R \iff x \in tR 
  \iff \hbox{the $A$-module $R/xA$ has torsion} 
  \iff \hbox{$R/xA$ is not flat} $$
  $$ x/t'\in R'\iff x \in t'R'
  \iff\hbox{the $A'$-module $R'/xA'$ has torsion}
  \iff \hbox{$R'/xA'$ is not flat.} $$

  So $x/t \notin R$ implies $R/xA$ is flat, which implies
  $R'/xA' = (R/xA) \otimes_{A}{A'}$ is flat, which implies $x/t' \notin R'$.

\junk{
  In the case at hand, $x/t\notin R$ so $xA$ is
  a cotorsion-free submodule of a finite-dimensional direct summand of $R$,
  the product of finitely many graded pieces. Since $A$ is a DVR, this
  implies that $xA$ is a direct summand of $R$. Then $xA'$ is a direct
  summand in $R'$, so $R'/xA'$ is torsion-free. Hence $x/t' \notin R'$.
  
  (The finite-dimensional assumption is not really necessary for this argument
  to work: $R/xA$ is torsion-free, so flat; hence, $R'/xA' = (R/xA)
  \otimes_{A}{A'}$ is flat, so torsion-free, since a module over a DVR is flat
  $\iff$ it is torsion-free $\iff$ it is free.)
}
\end{proof}

\begin{corollary}[The functor $\Bhb$ is separated]
  \label{cor:separated}
  An element of $\Bhb(K)$ has at most one extension
  to an element of $\Bhb(A)$.
\end{corollary}

\begin{proof}
  Call the original element $K[x_0,\dotsc,x_n]\to R_K$.
  The only possible $R$ will be the integral
  closure of the image of $A[x_0,\dotsc,x_n]$ in $R_K$,
  as we now show in two steps.

  Let $A[x_0,\dotsc,x_n]\to R$ be an extension. Then $R$ is finite
  over the image of $A[x_0,\dotsc,x_n]$ in $R$, so $\wt R$ is the integral
  closure of the image of $A[x_0,\dotsc,x_n]$ in $R_K$.
  
  The special fiber $R/tR$ is reduced because so is the geometric
  fiber $(R/tR)\otimes_k \bar k$. But then $R = \wt R$ by Lemma
  \ref{lem:nilpotents}(1). 
\end{proof}

\begin{remark}
  For any finitely generated algebra $S$ over a field $k$,
  $S\otimes_k \bar k$ is reduced $\implies$ $S$ is reduced. 
  The converse is true if $k$ is a field of characteristic zero or a
  perfect field of characteristic $p>0$.  Moreover, if $S\otimes_k \bar k$ 
  is not reduced then already for some finite purely inseparable
  extension $k'/k$, $S\otimes_k k'$ is not reduced.
\end{remark}

\begin{lemma}\label{lem:get-reduced} 
  Let $R$ be a Noetherian ring which has finite normalization (e.g.
  $R$ has no embedded primes and is finitely generated over a field
  or a DVR), $t\in R$ a nonzerodivisor, and assume that $R$ is
  integrally closed in $R[t\inv]$.  Then the ring $R/(t)$ does not
  have embedded primes; it satisfies Serre's condition $S1$.
  
  In particular, if $\Spec R$ is reduced, $R$ is finite over a DVR and
  integrally closed in $R[t\inv]$, and $\Spec R/(t)$ is generically
  regular, then $\Spec R/(t)$ too is reduced.
\end{lemma} 

By the normalization of $R$ we understand its integral closure in the total
ring of fractions $R[S^{-1}]$, the localization in all nonzerodivisors.

If we assumed $R$ to be reduced and equal to its normalization, then
$R$ would satisfy Serre's condition $S2$, so $R/(t)$ would be $S1$.  
This Lemma clarifies how much %of
integral closure one actually needs to draw this conclusion.

\begin{proof}[Proof of Lemma \ref{lem:get-reduced}, latter conclusion]
  There is a simpler proof in this special case
  (which is all we will need). Or rather, the only subtle point
  is encapsulated in a familiar formulation ``normal domains are $S2$''.
  
  Let $R'$ be the normalization of $R$. By our assumptions on $R$, the ring
  $R'$ is a direct sum of finitely many normal domains.  Consider the map
  $R/tR \to R'/tR'$. Its kernel is
  $$
  \frac{R \cap tR'}{tR} \cong
  \frac{t^{-1}R \cap R'}{R} \subset 
  \frac{R[t\inv]\cap R'}{R} = 0
  $$
  thanks to the assumptions on $t$ and on $R$.
  Hence the map $R/tR \to R'/tR'$ is an inclusion.
  
  Since $X_0 := \Spec R/(t)$ is generically regular, there exists an open
  subset $U\subset \Spec R$ such that $U$ is regular and such that
  $U\cap X_0$ is dense in $X_0$. Then the normalization morphism
  $\Spec R'\to \Spec R$ is an isomorphism over $U$. So $\Spec R'/(t)$
  is generically regular.

  Since $R'$ is $S2$, and $t$ is not a zero divisor, $R'/(t)$ is $S1$.
  Being $S1$ and generically regular, $R'/(t)$ is reduced.
  Since $R/(t)$ is a subring of a reduced ring, it too is reduced.
\end{proof}

\begin{proof}[Proof of Lemma \ref{lem:get-reduced} in general]
  Essentially, we must modify the proof that normal domains are $S2$.
  Let us introduce some notation.  Let $X=\Spec R$, and let $\mc F$ be
  a coherent sheaf on $X$, corresponding to a finite $R$-module
  $M$.  Let $X_0 = \Spec R/(t)$, and define the \defn{saturation of
    $\mc F$ in codimension 2 along $X_0$} by the formula
  \begin{displaymath}
    \mc F\sat = \varinjlim_Z \mc F(X\setminus Z),
  \end{displaymath}
  where $Z$ goes over closed subsets of $X_0$ that have
  $\codim_Z X_0\ge1$, equivalently $\codim_Z X\ge2$ since $t$ is a
  nonzerodivisor in $R$. Let $M\sat=\Gamma(X, \mc F\sat)$ be the
  corresponding $R$-module.

  It is true that the sheaf $\mc F\sat$ is coherent provided all $Z$s
  have codimension $\ge2$ in associated primes of $M$, cf.
  \cite[5.9-11]{EGA4-2}, but we do not need this. We merely observe that
  $\mc O_X\sat=\mc O_X$ since every $\Gamma(X\setminus Z,\mc O_X)$ is
  contained in the normalization of $ R$ and also in
  $\Gamma(X\setminus X_0, \mc O_X) =  R[t\inv]$.

  Let $G$ be a coherent sheaf supported on an irreducible subset $Z$ of $X_0$
  with $\codim_X Z\ge 2$. Then every extension
  \begin{displaymath}
    0\to \mc O_X \to F \to G \to 0
  \end{displaymath}
  splits. Indeed, $F\sat = \mc O_X\sat=\mc O_X$, and the canonical restriction
  morphism $F\to F\sat$ provides the splitting. Therefore, $\Ext^1(G,\mc
  O_X)=0$.  By the cohomological characterization of depth (see e.g.
  \cite[Thm. 28]{Matsumura_ComRingThry} or \cite[18.4]{Eisenbud_CommAlg}) this
  implies that the local ring $\mc O_{Z,X}$ has depth $\ge2$, and therefore
  $\mc O_{Z,X_0}$ has depth $\ge 1$, i.e. the ideal $p_Z$ is not an embedded
  prime of $R/(t)$.
\end{proof}

\begin{theorem}[The functor $\Bhb$ is proper]
  \label{thm:proper}
  Every element %$K[x_0,\dotsc,x_n]\to R_K$ 
  of $\Bhb(K)$ has an extension to one of $\Bhb(A')$, after a
  finite ramified base change $S'=\Spec A' \to S=\Spec A$.

  The necessary base change $\deg(S'/S)$ 
  divides $(\prod_{i=0}^{\deg h} b_i!)^2$,
  and if the residue field $k=A/(t)$ has characteristic zero, 
  then the necessary base change even divides $\prod_{i=0}^{\deg h} b_i!$.
  
  Let $X \to \Spec A$ be a flat proper extension (before any base
  change), $\wt X$ the normalization of $X$ in the generic fiber, and
  $\wt X_0$ the special fiber of $\wt X$.  Denote the multiplicities of
  the geometric fiber $\wt X_0\times_k \bar k$ by $\{m_i\}$.
  
  If $\chr k=0$, then the base change can be taken to be $t=s^m$, where
  $m=\lcm(\{m_i\})$, the least common multiple of the multiplicities.  

  If $\chr k>0$, the base change can be chosen to be a composition of a base
  change $A'/A$ such that $tA'=t'A'$ and the residue field extension
  $k'/k$ is purely inseparable of degree dividing $\prod m_i$, and the
  base change $t'=s^m$.
\end{theorem}

\begin{proof}
  We first provide a flat $A$-algebra $R$ extending $R_K$.  Let $q_1,\dotsc,
  q_s\in R_K$ be homogeneous elements generating $R_K$ as a
  $K[x_0,\dotsc,x_n]$-module. Each $q_i$ is integral over $K[x_0,\dotsc,x_n]$.
  Since $K=A[1/t]$, there exist $n_i\in\naturals$ such that $t^{n_i}q_i$ are
  integral over $A[x_0,\dotsc,x_n]$. Therefore, the algebra $R_1:=(\im
  A[x_0,\dotsc,x_n])[t^{n_i}q_i]$ is finite over $A[x_0,\dotsc,x_n]$,
  graded, and free over $A$. 
  Take $R=\wt R_1$ to be its normalization in $R_1[t\inv]$, as before.
  (These $R_1,\wt R_1$ are the coordinate rings of $X,\wt X$ in the
  statement of the Theorem.)

  Next, we find a ramified base change $S'\to S$ so that the
  geometric fiber $\Spec \wt{R'}\otimes_{k'}\bar k'$ is reduced. By
  Lemma~\ref{lem:get-reduced} we only need to prove that it is generically
  reduced.

  Assume first that $\chr k=0$.
  Let $Z$ be an irreducible component of the central fiber $X_0$, with 
  generic point $z$. The ring $\mc O_{z,X}$ has dimension $1$ and is integrally
  closed in the generic fiber, since normalization commutes with localization.
  (This is where we must use $\wt R_1$ and not $R_1$.)
  Therefore, it is integrally closed, so it is a DVR with a uniformizing
  parameter, denote it by $\pi$.
  
  We have $t=a\pi^n$ for some invertible $a$. Let us make a base
  change $t=s^n$. Then the normalization $\widetilde{\mc O}_{z',X'}$
  of $\mc O_{z,X}\otimes_A A'$ is regular and its central fiber has
  multiplicity~1.  Indeed, $\pi/s \in \widetilde{\mc O}_{z,X}$ and is
  invertible, so $s= b\pi$ in $\widetilde{\mc O}_{z',X'}$ with
  invertible $b$.  Hence, after making the base change $t=s^m$, where
  $m=\lcm(\{m_i\})$, the central fiber is generically regular, and we
  are done.
  
  If $\chr k>0$, then we first make an unramified base change $S'/S$ with a
  purely inseparable field extension $k'/k$ after which the multiplicities of
  the irreducible components of $X_0'$ and $X_0\times_k \bar k$ are the same;
  the degree of this base change divides $\prod m_i$. Then we proceed as
  above. 
\end{proof}

\junk{
\begin{proof}[Second proof.]
  If the residue field $k$ is big enough in the sense of
  Remark~\ref{rem:big-field} (for example, $k$ is algebraically closed) and
  has characteristic $0$, there is an alternative, more geometric proof. 
  
  In this case, for each dimension $i$ we have a reduced $1$-dimensional scheme
  $S_i=(X_i)^{\dim 0}$ with a finite morphism to $S$.  Since $R$ is integrally
  closed in the generic fiber $R[1/t]$, the generic linear section $S_i$ is
  also integrally closed in its generic fiber.  Hence, $S_i$ is a disjoint
  union $\coprod S_{ij}$ of $1$-dimensional regular schemes, spectra of DVRs,
  with finite morphisms $f_{ij}:S_{ij}\to S$. % of degrees $b_{ij}$, $\sum_j
%  b_{ij}=b_i$.
  
  After making the base change $S'=S_{ij}\to S$, $S'\times_S S'$ acquires a
  section. Hence, after making finitely many base changes $S_{ij}\to S$ and
  normalizing in the generic fiber we get a section through every irreducible
  component of the central fiber. Consequently, the geometric central
  fiber becomes generically regular. The product of the degrees of the
  base changes obviously divides $\prod_i b_i!$.
\end{proof}
}

\begin{corollary}
  For any further ramified base change $A''\supset A'$, the
  central fiber does not change except for the extension of the
  residue fields:  $X''_0 = X'_0 \otimes_{k'} k''$.
  In other words, every 1-parameter family of branchvarieties has a unique
  limit, up to extensions of the residue field.
\end{corollary}

\begin{proof}
  Indeed, the geometric central fiber of $X' \otimes_{S'} S''$ is the same as
  that of $X'$, so it is reduced. By Lemma~\ref{lem:nilpotents}(1), one has
  $X'' = X'\times_{S'} S''$.
\end{proof}

The $\{m_i\}$ have a simple interpretation in the case that $X$ is a
family of points over $\Spec \mb C[[t]]$: 
they are the lengths of the cycles in the monodromy around $0$ 
of the generic fiber.
The degree $\lcm(\{m_i\})$ base change replaces the monodromy in this
finite family by a high enough power to make it trivial.

\begin{example}
  In examples, where one has a family of subvarieties limiting to a
  subscheme, one often knows the Hilbert extension $X\to \Spec A$.
  So it is tempting to work with the multiplicities of the special
  fiber $X_0$ to compute the necessary base change, rather than those of
  $\wt X_0$ as stated in the Theorem.  But consider the family 
  $$X = \left\{[x,y] : (y^2 + t x^2)(y^3 - t x^3) = 0\right\}$$ of
  $5$-tuples of points.  In this case $X_0$ is a quintuple point, but
  the necessary base change has $\deg(S'/S)=6$, correctly calculable
  from $\wt X_0 = \{$double point, triple point$\}$.
\end{example}

\begin{remark}
  In \cite{Alexeev_Limits} it is proved that the canonical limits of varieties
  and pairs of general type, whose construction follows from the log Minimal
  Model, are $S2$, in addition to being $S1$, as in Theorem~\ref{thm:proper}.
  This result, however, applies to a less general class of degeneration
  families than \eqref{thm:proper}.
\end{remark}

%%%%%%%%%%%%%%%%%%%%%%%%%%%%%%%%%%%%%%%%
\section{Construction of the moduli space}
\label{sec:constr-moduli-space}

For every family of branchvarieties $f:X\to Y_S$, define a functor
$$\Aut(f): (S\text{-schemes}) \to (\text{Groups})^{\op} $$
by setting $\Aut(f)(S')$ to be the automorphism group of $f':X'\to Y'_{S'}$,
where $X'=X\times_S S'$, $Y'=Y\times_S S'$, and $Y'_{S'}=Y_S \times_S S' =
Y\times_C S'$.

\begin{theorem}\label{thm:finite-aut}
  $\Aut(f)$ is represented by a finite group scheme over $S$.
\end{theorem}
\begin{proof}
  Let $X^{(2)}:=X\times_{Y_{S}} X$, a proper and projective scheme over $S$.
  
  For every automorphism $g:X'\to X'$ over $Y'_{S'}$, its graph $\Gamma_g$ is
  a closed subscheme of $X^{(2)}\times_S S'$. Therefore, it represents an
  $S'$-point of the Hilbert scheme $\Hilb X^{(2)}$, i.e.  an element of
  $(\Hilb (X)^2)(S')$. Moreover, there is a natural open subscheme $U$ of
  $\Hilb X^{(2)}$ parametrizing subschemes $Z\subset X^{(2)}\times_S S'$ that
  project isomorphically to both factors, and $\Gamma_g$ gives an $S'$-point
  of $U$. The opposite is clear as well: every such subscheme $Z$ is the graph
  of a unique automorphism. Hence, the quasiprojective-over-$S$ scheme $U$
  represents $\Aut(f)$. It is obviously a group scheme.
  
  $\Aut(f)$ satisfies the valuative criterion of properness thanks to
  the properness of the functor of branchvarieties,
  Theorem~\ref{thm:proper}.  To prove that $\Aut(f)$ is finite over
  $S$ we need to check that the geometric fibers are finite.
  
  So let $f:X\to Y$ be a branchvariety over a field $k$. Cover $Y$ by
  finitely many open affines $V_i$. Since the morphism $f$ is finite,
  each $U_i= f\inv(V_i)$ are affine as well, and we only need to show
  that $\Aut(U_i/V_i)$ is finite. Let $r\in k[U_i]$. It satisfies some
  monic polynomial equation with coefficients in $k[V_i]$, and the
  image of $r$ under any automorphism must satisfy it too. Since
  $k[U_i]$ is reduced, it is embedded into a direct sum of finitely
  many fields, one for each irreducible component. A monic (hence nonzero)
  polynomial has only finitely many roots in a field, so we are done.
\end{proof}

For a branchvariety $X$ over a field of characteristic $0$ it is easy to say
more: $\Aut(f)$ is a subgroup of the product of Galois groups of the
irreducible components of $X$, and therefore a subgroup of the product of
several symmetric groups. And, of course, by Cartier's theorem 
any group scheme in characteristic zero is reduced.

\begin{example}\label{ex:nonreduced}
  Let $\chr k=p$ and $f:X=\bP^1\to Y=\bP^1$ be the geometric Frobenius
  morphism, $x\mapsto x^p$. Then $\Aut(f)=\mu_p=\Spec k[x]/(x^p-1) = \Spec
  k[x]/(x-1)^p$, a finite nonreduced group scheme.
\end{example}

\medskip
\begin{proof}[Proof of the Main Theorem~\ref{thm:main}]
  We first prove that branchvarieties together with some additional data can
  be parametrized by a locally closed subscheme of a certain Hilbert
  scheme. This is done by a classical argument, as in the case of curves 
  \cite[Prop. 5.1]{Mumford_GIT}, with necessary modifications. 
  
  Then the moduli space is constructed by taking the quotient by a $\PGL$
  group action. We do not use Geometric Invariant Theory for this step.
  Instead, the quotient by a proper group action immediately gives the moduli
  stack as an algebraic Artin stack. By applying standard results on
  representability, we obtain its coarse moduli space as an algebraic space.
  
  Let $f:X\to Y$ be a branchvariety defined over an algebraically closed field
  $k$, with fixed Hilbert polynomial $h$ and degree sequence $b$.  By Theorem
  \ref{thm:boundedness}, we know that there exists some integer $d_0(h,b)$
  such that $L^{d_0}$ has vanishing higher cohomology, is very ample, and such
  that the ring of global sections $R(X,L^{d_0})$ is generated in degree 1.
  Let $D=h^0(X,L^{d_0})-1$.  A choice of a basis in the vector space
  $H^0(X,L^{d_0})$ defines an embedding $X\into \bP^D$, and two such choices
  differ by an element of $\PGL_{D+1}(k)$. 
  
  Let $\bP^n\supset Y$ be the projective embedding of $Y$.
  Let $\bP^D\times \bP^n\into \bP^m$ be the Segre embedding, 
  $m+1 =(D+1)(n+1)$. The restriction of $\mc O_{\bP^m}(1)$ to 
  $X\subset \bP^D\times Y \subset \bP^D\times \bP^n \subset \bP^m$ is 
  isomorphic to $L^{d_0}\otimes L$ and has Hilbert polynomial $H(d) =
  h(d(d_0+1))$.  
  
  Let $\Hilb_{H,\bP^D\times Y}$ be the Hilbert scheme parametrizing closed
  subschemes of $\bP^D\times Y$ that have Hilbert polynomial $H$.  The
  properties of being geometrically reduced, the first projection being a
  closed embedding with the image spanning $\bP^D$, and the second projection
  being finite, are all open in projective families over a quasiprojective
  base.  Thus, there exists an open subscheme $V_1$ of $\Hilb_{H,\bP^D\times
    Y}$ whose $k$-points correspond to branchvarieties of $Y$ embedded in
  $\bP^D\times Y$ and spanning $\bP^D$.
  
  The Hilbert polynomials of the ample sheaves $p_1^*\mc O_{\bP^D}(1)$ and
  $p_2^*\mc O_{Y}(1)$ are locally constant.  By Lemma
  \ref{lem:degree-sequence}, the degree sequence is locally constant as well.
  This gives an open subscheme $V_2 \subseteq V_1$ over which the 
  branchvarieties have invariants $(h,b)$ and such that the sheaf
  $\cO_{\bP^D}(1)|_X$ has the same Hilbert polynomial as the sheaf $L^{d_0}$.
  
  Let $\mc X_2\to V_2$ be the universal family. On this family, we have
  two invertible sheaves, $p_1^*\mc O_{\bP^D}(1)$ and $p_2^*\mc O_{Y}(d_0)$.
  We claim that there exists a locally closed subscheme $V$ of $V_2$
  parametrizing branchvarieties on which these two sheaves coincide. Indeed,
  the relative Picard functor $\Pic_{\mc X_2/V_2}$ is represented by an
  algebraic space (this is Artin's theorem, see
  \cite[Thm. 7.3]{Artin_FormalModuli1} or
  \cite[Thm. 8.3.1]{BoschLutkebohmertRaynaud90}). The two sheaves above define
  sections of this algebraic space over $V_2$, and $V$ is the locus where
  these sections coincide.
  
  To summarize: we have constructed a scheme $V$ parametrizing
  branchvarieties of $Y$ together with an embedding by a complete linear
  system $L^{d_0}$ into a fixed projective space $\bP^D$ so that the image
  spans $\bP^D$.  Two points in $V(k)$ define isomorphic branchvarieties iff
  they are in the same orbit of the group $\PGL_{D+1}(k)$.
  
  Now consider a family of branchvarieties $f:X\to Y_S$, $\pi:X\to S$ over an
  arbitrary $C$-scheme $S$. By the Cohomology and Base Change Theorem
  \cite[III.12.11]{Hartshorne}, 
  the sheaf $F_{d_0}=\pi_* L^{d_0}$ is locally free,
  so it becomes trivial on some open affine cover $S=\cup S_i$. The choice of
  trivializations of $F_{d_0,i}=F_{d_0}|_{S_i}$ gives the choice of embeddings
  of $X_i$ into $\mb P^D \times Y_{S_i}$, and two such embeddings differ by an
  element of $\PGL_{D+1}(S_i)$. This gives a collection of $S_i$-points of
  $V$, up to actions of the groups $\PGL_{D+1}(S_i)$.
  
  It follows that the stack in groupoids $\cB^b_{h,Y}$ is just the quotient
  stack $[V / \PGL_{D+1} ]$, in other words the quotient of $V$ by a smooth
  pre-equivalence relation $$j:R=V\times \PGL_{D+1}\to V\times V.$$ 

  It is well known that the separatedness of the moduli functor (Corollary
  \ref{cor:separated}) and finiteness of the automorphism schemes
  (Theorem~\ref{thm:finite-aut}) imply that the group action is proper. 
  In particular, the stabilizer  $j\inv(\diag V)\to V$ is finite.
 
  Thus, $\Bhb$ is an algebraic Artin stack; see \cite{LaumonMoretMailly} for a
  general reference. By either \cite{KeelMori} or \cite{Kollar_QuotSpaces} it
  has a coarse moduli space as an algebraic space (see \cite{Knutson_AS} for
  the general reference on algebraic spaces).  Finally, by
  Theorem~\ref{thm:proper} the functor is proper, so the moduli space is too.
\end{proof}

\begin{remark}
  A separated Artin stack with a finite stabilizer is
  \defn{Deligne-Mumford} if its stabilizer groups are \emph{reduced}.
  By Theorem \ref{thm:finite-aut}, in characteristic $0$ our moduli
  stacks of branchvarieties are Deligne-Mumford (and more generally,
  under Assumption \ref{assume:char} to come), but Example
  \ref{ex:nonreduced} shows that in general they are not.
  Nor are the stabilizer groups always linearly reductive in
  characteristic $p$, as the branchvariety of $p$ points mapping to a point
  (with automorphism group $S_p$) demonstrates.
\end{remark}

\begin{remark}
  Recall that the tangent space to the Hilbert scheme has a particularly
  simple description: if $Z\subset {\mb P}^D \times Y$ 
  is a closed subscheme defined by an ideal sheaf $I_Z$ then 
  $$ T_{[Z],\,\Hilb({\mb P}^D \times Y )} = \Hom( I_Z/I_Z^2, \mc O_Z). $$

  Applying \cite[Theorem 1.5]{Olsson} to the representable morphism
  $V \onto [V/\PGL_{D+1}]$, and chasing a couple of short exact sequences
  relating 
  the cotangent complexes \cite{Illusie_Cotangent1,Illusie_Cotangent2},
  we can describe the tangent space to the corresponding point of $\B$ as
  $$ T_{[Z],\,\B(Y)} 
  = T_{[Z],\,\Hilb(\mb P^D \times Y)} \,\big/\, Lie(\PGL_{D+1}). $$
  In addition, one can identify the corresponding obstruction spaces
  on the nose (no quotient by $Lie(\PGL_{D+1})$).

\junk{  
  An infinitesimal deformation of a finite morphism is obviously
  finite. Hence, the tangent space to $\Bhb$ can be computed by using the
  general machinery of cotangent complexes
  \cite{Illusie_Cotangent1,Illusie_Cotangent2}. The additional condition that
  $f:X\to Y$ is finite does not seem to simplify the situation much.

  On the other hand, our construction of the moduli stack as the quotient
  $[V/\PGL_{D+1}]$ gives a description of the tangent space in terms of
  the tangent space to a Hilbert scheme and the Lie algebra of $\PGL_{D+1}$.
}
\end{remark}

%%%%%%%%%%%%%%%%%%%%%%%%%%%%%%%%%%%%%%%%%%%%%
\section{Examples}
\label{sec:examples}

We begin with some simple examples in the first section, and then present
the ones that principally guided each of us to the theory of branchvarieties.

%%%%%%%%%%%%%%%%%%%%
\subsection{Further elementary examples}

In the examples below, we take $A= \mb C[[t]]$, with $K$ its field of
fractions. 

\begin{example}\label{ex:skewlines}
  Consider two skew lines in $\AA^3$ approaching each other as $t\to 0$:
  $R = A [x,y,z] / (z,x)\cap (z-t,y)$. The central fiber
  $$  A[x,y,z] / (xy, z^2, zx, zy) $$
  is two lines with an embedded prime
  at the point of intersection, and is not reduced.  
  No ramification is necessary in this case (as follows 
  from Theorem \ref{thm:proper}, since the central fiber is
  generically reduced).
  The integral closure is the union of two disjoint families of lines
  $$  \wt R= A [x,y,z]/(z,x) \oplus A [x,y,z]/(z-t,y).$$
  The central fiber is a disjoint union of two lines, with a finite map to the
  two intersecting lines in $\AA^3$.
  
  The total space of the Hilbert family $\Spec R$ is the union of two planes
  meeting at a point, everyone's first example of a scheme smooth in
  codimension $1$ yet still abnormal: it is not $S2$, and a hyperplane
  section has an embedded prime.  The integral closure just pulls
  those two planes apart.
\end{example}

\newcommand\Pone{{\mb P}^1}

\begin{example}\label{ex:3lines}
  Let $X$ be a union of 3 copunctal lines in $\mb P^2$, say with slopes 
  $0$, $1$, and $\infty$, and let $f_t:X\to \mb P^1$ be a projection
  with the angle $t$. The tangent space to the common point is
  $2$-dimensional, and contains four $1$-d subspaces: the tangents to
  the lines, and the kernel of the derivative of $f_t$. In these terms,
  $t$ is a cross-ratio. Consider this as a 1-parameter family as $t$
  goes to $0$. 
  
  In coordinates, the generic fiber $X_{\eta}$ corresponds to a graded ring
  $R$ which is a subring of $K[x_0,x_1] \oplus K[y_0,y_1] \oplus K[z_0,z_1]$
  consisting of homogeneous polynomial $(f,g,h)$ such that
  \begin{displaymath}
    f(1,0)=g(1,0)=h(1,0), \qquad tf'(1,0)+g'(1,0)= (1+t)h'(1,0),
  \end{displaymath}
  where the derivatives are with respect to the second variable, $x_1$, $y_1$
  or $z_1$ respectively. The $K[t_0,t_1]$-module structure is given by the
  homomorphisms $t_i\mapsto (x_i,y_i,z_i)$. 
  
  We can use the same equations for the total family $X$, replacing $K$ by
  $A$.  Specializing $t=0$, one obtains the ring of triples $(f,g,h)$ such
  that
  \begin{displaymath}
    f(1,0)=g(1,0)=h(1,0), \qquad  g'(1,0)= h'(1,0),
  \end{displaymath}
  Hence, the central fiber $X_0$ is a union of three $\mb P^1$s passing
  through one point, the second and the third $\mb P^1$s are tangent to each
  other, and the first is transverse to them. Note that $X_0$ can no longer be
  embedded in $\mb P^2$.
\end{example}

\begin{example}\label{ex:complete-conics}
We will assume for simplicity that $\chr k \neq 2$ in this example.

Consider the space of branchvarieties of $\PP^2$ including the plane
conics, so $h(n) = 2n+1$, $b_0 = 0$, $b_1 = 2$. In short, this
coarse moduli space is the ``space of complete conics'', or $\PP^5$
blown up along the Veronese surface, while the stack structure agrees
with that on the corresponding Kontsevich moduli space of stable maps.
Indeed, each branchcurve $X$ has arithmetic genus 0, and so has at worst nodes
as singularities.
We include some standard facts about this stack.

There are two obvious closed substacks: $T = \{$reducible branchvarieties$\}$
and $N=\{$noninjective branchvarieties$\}$. Each of the points in the substack 
$N$ has automorphism group ${\mathbb Z}/(2)$. The complementary set $N^c$
corresponds to reduced plane conics, has trivial stack structure, 
and is isomorphic as a space to $\PP^5$ with the Veronese surface removed.
The set $T^c \cap N$ consists of double covers $\PP^1 \to \PP^1$
of lines; such a cover is uniquely determined by its image line and the
two (distinct) branch points. If we let the branch points collide, we
get the space $T\cap N$ consisting of pairs of crossing lines double
covering a line in $\PP^2$. 

The whole space is $5$-dimensional. The substack $N$ is $4$-dimensional.
The substack $T\cap N$ is $3$-dimensional, and isomorphic as a stack to the
manifold of flags in $\PP^2$ modulo a trivial ${\mathbb Z}/(2)$ action.

In the case of plane cubics, the moduli stacks of branchcurves and
of stable maps are \emph{not} naturally isomorphic, which we will see
in Section \ref{sec:complete-cubics}.
\end{example}

\begin{example}\label{ex:double-covers}
Fix $g\in\naturals$, and let $h=2n+1-g$, $b_0=0$, $b_1=2$.
We can make some of the branchcurves of $\PP^1$ with these invariants 
by joining two $\PP^1$s along $g+1$ distinct points; the arithmetic
genus of such a curve is $g$.

The generic branchcurve of $\PP^1$ with these invariants is a smooth curve 
of genus $g$, branched over $\PP^1$ at $2(g+1)$ points, and the coarse 
moduli space is just $$(\PP^1)^{2(g+1)}/S_{2(g+1)} \cong \PP^{2(g+1)}.$$

For example, consider $g=1$, so the generic branchcurve is an elliptic
curve branched over $\PP^1$ at four distinct points. If two branch points 
coalesce, $4 = 2+1+1$, the branchcurve is a nodal cubic with the node mapping
to the double branch point. If the other two coalesce also,
$4 = 2+2$, the branchcurve is as described a moment ago -- two $\PP^1$s
glued together at two points. If three branch points coalesce,
$4 = 3+1$, the branchcurve is a cuspidal cubic with the cusp
mapping to the triple branch point. If all four coalesce, 
the branchcurve is a union of two $\PP^1$s along a point of tangency.

We note that the moduli space of branchvarieties of $\bP^1$ contains 
the classical Hurwitz schemes parametrizing degree $d$ covers of $\bP^1$
with certain ramification conditions. We see from the above example that the
compactification of these Hurwitz schemes provided by branchvarieties is very
different from the compactification obtained by adding ``admissible covers''
as done, e.g., in \cite{HarrisMumford_KodairaDim}.
\end{example}

%%%%%%%%%%%%%%%%%%%%
\subsection{Stable toric varieties}
\label{sec:stable-toric-vari}

Let $T=(\mb G_m)^r$ be a split multiplicative torus, a direct sum of $r$ copies
of the multiplicative group $\mb G_m$. Let $\bP^n$ be a projective space
endowed with a $T$-action and with a $T$-linearized $\mc O(1)$.

\begin{definition}
  A \defn{stable toric variety over $\bP^n$} over an algebraically closed
  field is a seminormal projective variety $X$ 
  endowed with a $T$-action such that
  \begin{enumerate}
  \item there are only finitely many orbits, and
  \item the isotropy groups are subtori, so in particular connected and
    reduced
  \end{enumerate}
  together with a finite and $T$-equivariant morphism $f:X\to \bP^n$.
\end{definition}

When $T$ acts on $(\bP^n,\mc O(1))$ with $n+1$ distinct characters, the data
for the morphism $f:X\to \bP^n$ is equivalent to the data for an effective
ample Cartier divisor $D$ on $X$ which does not contain any $T$-orbits. So, in
this case a stable toric variety over $\bP^n$ is the same as a stable toric
pair $(X,D)$ with a $T$-linearized line bundle $L=\mc O_X(D)$, see
\cite[2.14]{Alexeev02}, \cite[Prop. 3.3.2]{AlexeevBrion05}. The latter paper
also includes the more general case of stable spherical varieties, where $T$ 
is replaced by a reductive group.

The higher cohomologies $H^i(X,L^d)$ vanish for $d>0$. Therefore, for every
family $X\to \bP^n_S$ of stable toric varieties the graded $\mc O_S$-algebra
$ R(X/S, L) = \oplus_{d\ge 0} \pi_* L^d$
is locally free. 
The $T$-action is equivalent to the grading of $R(X/S,L)$ by
the character group $\Lambda=\mb Z^r$ of $T$, and each graded piece
$R(X/S,L)_{\lambda}$, $\lambda\in \Lambda$, is of finite rank $h(\lambda)$,
i.e. $R(X/S,L)$ is multiplicity-finite. 
When $R(X/S,L)$ is \defn{multiplicity-free}, i.e. each $h(\lambda)$ is $0$ or
$1$, a stronger statement is true: $X$ reduced implies that $X$ is
seminormal. 

We see that the moduli of multiplicity-free stable toric varieties over
$\bP^n$ is just the branch-analogue of the toric Hilbert scheme of
Peeva-Stillman \cite{PeevaStillman} and the more general multigraded Hilbert
scheme of Haiman-Sturmfels \cite{HaimanSturmfels}. 

This ``toric Branch'' moduli space is constructed in \cite{AlexeevBrion05}
more directly, as the quotient $U/\Gamma$ 
of $U= \Hilb_h(Z)$, where $Z\to \bP^n$ is an $A$-cover,
by a finite diagonalizable group,
i.e. product of several groups $\mu_m$ of roots of unity. 
The reason for the relative simplicity of this case is that
the monic polynomials appearing in the finite ring extensions have the very
simple form $z^m=r$. The projectivity of the moduli space is
immediate from this description.

Since some of the varieties appearing below in Section \ref{sec:u-invts} are
stable toric varieties, we recall briefly their classification.  Each stable
toric variety over an algebraically closed field defines a complex of
polytopes $\Delta$ with a reference map to $\Lambda_{\mb R}$.  This means that
we have a topological space $|\Delta|$ with a cell decomposition $|\Delta|
=\cup \delta$ and a finite map $\rho:|\Delta|\to \Lambda\otimes\mb R$
identifying each $\delta$ with a lattice polytope. Then $X$ is a union of
ordinary (normal) projective toric varieties $X_{\delta}$ which are glued the
same way as the complex $\Delta$.

A variety $X$ is multiplicity-free precisely when the map $\rho$ is
injective. One-parameter degenerations correspond to convex subdivisions of
$\Delta$.

%%%%%%%%%%%%%%%%%%%%
\subsection{Balanced normal cones}
\label{sec:balanc-norm-cones}

\newcommand\gr{{\mathop {\mathrm gr}}}
\newcommand\ogr{{\mathop {\overline{ \mathrm gr}}}}

Let $Q$ be a commutative $k$-algebra without nilpotents, and $I$ an ideal.
The \defn{Rees algebra} is the graded subring 
$$ R = \left(\bigoplus_{n<0} t^{-n} I^n\right) 
\oplus \left( \bigoplus_{n\geq 0} t^n Q\right) $$
of $Q[t,t^{-1}]$. Under the evident map $k[t] \to R$, we see that the
map $\Spec R \to \Spec k[t]$ defines a flat family whose $t=1$ fiber is $Q$ 
and $t=0$ fiber is $\gr_I Q := Q/I \oplus I/I^2 \oplus I^2/I^3 \oplus \ldots$, 
the associated graded with respect to the $I$-adic filtration. If we
assume in addition that $\cap_n I^n = \{0\}$, then this family is
locally free.  Geometrically, this family is the 
\defn{degeneration of $\Spec Q$ to the normal cone} of $\Spec Q/I$. 
% in $\Spec Q$.

\newcommand\barf{{\overline f}} In \cite{Samuel}, Samuel defined a
variant of the $I$-adic filtration
$$ \forall q\in Q,\quad  f(q) := \max \{n : q \in I^n \} $$
called its \defn{homogenization},
$$ \forall q\in Q,\quad  \barf(q) := \lim_{k\to\infty} \frac{f(q^k)}{k}, $$
and proved that this limit exists. Rees (see the book \cite{Rees1})
and Nagata \cite{Nagata} proved that the limit is rational with
bounded denominator (depending on $Q,I$).  Let $N$ be divisible by all
the possible denominators of $\barf$; of course their LCM will do.

To this ``homogeneous'' filtration $\barf$, one can again associate a
Rees algebra (now $\frac{1}{N} \naturals$-graded) giving a flat
degeneration of $Q$ to an associated graded ring $\ogr_I R$,
this time \emph{automatically without nilpotents}. 
The corresponding geometry was studied in \cite{Knutson05} under the
name ``degeneration to the \emph{balanced} normal cone''.

We now relate this construction to the one in Section
\ref{sec:one-param-families}. 
Make the ramified base change $t = (t')^N$. 
Let $\wt R$ denote the integral closure of $R$ in $R \tensor_{k[t]} k[t']
\subseteq Q[t,t^{-1}]\tensor_{k[t]} k[t'] = Q[t',t'^{\, -1}]$. 
Then for $n\geq 0$,
\begin{eqnarray*}
  q/t'^{\, n} \in \wt R 
  &\Longleftrightarrow& q/t'^{\, n} \hbox{ is integral over $R$ }\\
  &\Longleftrightarrow& (q/t'^{\, n})^N \hbox{ is integral over $R$ }\\
  &\Longleftrightarrow& q^N/t^n \hbox{ is integral over $R$ }\\
  &\Longrightarrow& \barf(q) \geq n/N.
\end{eqnarray*}
By Rees' valuative formula for $\barf$ \cite[Thm. 4.16]{Rees1},
the converse of this last implication is also true. Hence
$$ \wt R = \left( \bigoplus_{n<0} t'^{\, -n} \{q : \barf(q)\geq n/N \}\right)
\oplus \left( \bigoplus_{n\geq 0} t'^{\, n} Q\right) $$
is the Rees algebra associated to the filtration given by $\barf$.

Craig Huneke informed us that much the same interpretation of Rees'
results on the Samuel filtration occurs in Theorem 10.6.6 of his
forthcoming book \cite{HuSw}.

\subsection{Chiriv\`\i's degeneration of flag manifolds as 
  a limit of branchvarieties}
\label{sec:Chirivi}

We describe a special case of Chiriv\`\i's geometric interpretation of
the Littelmann-Lakshmibai-Seshadri weight multiplicity formula \cite{Chirivi}, 
using the language of balanced normal cones. This example was what
motivated the second author to seek a general theory of automatically
reduced degenerations. The details will appear elsewhere \cite{Knutson06}.

Let $G$ be a complex connected algebraic group with maximal torus $T$,
and $\lambda$ a dominant weight. Then there is a natural $G$-equivariant
graded ring structure on $R := \oplus_{n\in\naturals} V_{n\lambda}$
(the $n$th piece being the irreducible representation of $G$ with high
weight $n\lambda$), whose $\Proj$ is a generalized flag manifold $G/P$.

Through a careful analysis of generators and relations of the ring $R$,
Chiriv\`\i\ gave a collection of $T$-equivariant degenerations $R'$ of $R$, 
where each $\Proj R'$ is a stable toric variety.  In some cases,
the underlying complex $\Delta$ of polytopes is in fact a simplicial
complex, with one simplex for each chain in the Bruhat order of $G/P$;
we will call this a {\em simplicial} Chiriv\`\i\ degeneration.

This degeneration was already well known in the case that $G/P$ is a
Grassmannian in its Pl\"ucker embedding \cite{DEP}. In this case,
each component of the stable toric variety maps isomorphically, not
just finitely, to a coordinate subspace of projective space.

By the flatness of these degenerations, one can compute the $T$-weight
multiplicities in the representation $V_\lambda$ as a sum over chains
in the Bruhat order, and for each chain, a count of lattice points in
a certain simplex (determined using $\lambda$). This weight multiplicity
formula had already been proven by Littelmann using his path model,
confirming a conjecture of Lakshmibai and Seshadri (inspired by
\cite{DEP} and followup work \cite{DL} generalizing it to other
minimal embeddings of minimal classical flag manifolds).

We now sketch a way to see a simplicial Chiriv\`\i\ degeneration as a
flat limit of branchvarieties, with proofs to appear in \cite{Knutson06}.  
The principal benefit of this viewpoint is that the construction does
not require special analysis of the ring $R$.

\newcommand{\Sym}{\operatorname{Sym}}

The {\em extremal weights} of $V_\lambda$ are of the form $w\cdot\lambda$ for 
$w$ in the Weyl group of $G$. Each extremal weight space is $1$-dimensional; 
let $E\leq V_\lambda$ be their direct sum. 
Then we make $G/P = \Proj R$ a branchvariety of projective space,
using the map $\Proj R \to \Proj \Sym E$.
(Indeed, $E$ is the smallest $T$-invariant subspace such that this map
has no basepoints.)

\newcommand\degen{\ \ \ \dashrightarrow\ \ \ }

The Bruhat order of $G/P$ gives a natural partial order on the
extremal weights; pick a linear extension (which will be immaterial)
to a total order.
Running through the sequence of extremal weights, 
we get a series of degenerations to balanced normal cones
$$ R \degen \ogr_{I_1} R \degen \ogr_{I_2} \ogr_{I_1} R  
\degen \ldots \degen \ogr_{I_m}\!\! \cdots\ogr_{I_1} R  $$
where the ideal $I_k$ is generated by the $k$th extremal weight space in the
total order (or rather, by the image of that extremal weight space 
in the $(k-1)$st ring in the sequence).

Since these ideals are $T$-invariant, the degenerations are $T$-equivariant.  
Each component of the resulting scheme is a weighted cone on a
weighted cone on $\ldots$ on a point, i.e. a toric variety associated
to a weighted simplex.

It is reasonably straightforward to show that the final ring so
constructed is a subring of Chiriv\`\i's ``discrete LS algebra'' $R'$.
One can then e.g. invoke Littelmann's result to show they are equal.

%%%%%%%%%%%%%%%%%%%%

%%%%%%%%%%%%%%%%%%%%%%%%%%%%%%%%%%%%%%%%
\section{Line bundles on Branch}
\label{sec:line-bundles-moduli}

By Section \ref{sec:boundedness}, for any $d\ge d_0$ the sheaves $L^d$ on our
branchvarieties do not have higher cohomology. Therefore, for any family
$f:X\to Y_S$, $\pi:X\to S$, the sheaves $F_d= \pi_*L^d$ are locally
free, and they induce natural line bundles $\lambda_d = \det F_d$ on our
moduli stacks. 
One might therefore hope (cf. \cite{Kollar_Projectivity}
for a quite similar situation) that these line bundles are ample. %for $d\gg0$.
(In particular, this would immediately imply that the coarse moduli
spaces of branchvarieties are projective schemes.) So, the following
is a somewhat surprising observation.

\begin{example}
  Consider the moduli space of 1-dimensional branchvarieties of $\bP^2$
  which include the plane cubics. For any family $f:X\to \bP^2_S$ whose fibers
  are reduced planar cubics, each $\lambda_d|_S$ has a positive
  degree. Indeed, this part of the branchvariety moduli coincides with 
  part of the Hilbert scheme, and $\lambda_d$ is the pullback of the
  (very ample) standard line bundle on a Grassmannian into which the
  Hilbert scheme is embedded.
  
  On the other hand, consider the family over $S=\mb P^1$ of
  Example~\ref{ex:3lines} (see also Section~\ref{sec:complete-cubics}). In
  this family, $\mc O_X$ is embedded into $\mc O_{\wt X}=\oplus_{i=1}^3 \mc
  O_{\mb P^1}$, for the normalization $\wt X$ of $X$, and the induced morphism
  $\wt X\to Y=\bP^1$ is constant.  Therefore, each $F_k$ is a nonconstant
  subbundle of a constant vector bundle 
  $\oplus_{i=1}^3 H^0(\bP^1,\mc O(d))\otimes \mc O_{S}$.  
  Hence, $\lambda_d|_S$ has a negative degree!
\end{example}

\begin{remark}
  We \emph{do not} claim that the moduli spaces $\Bhb$ are not projective. In
  particular, it seems possible that for $d,k\gg 0$ the line bundles 
  $\lambda_{dk}\otimes{\lambda_d}^{-(\deg h-1/2)k}$ on $\Bhb$ are ample.
  See also Theorem \ref{thm:Kollar}, which gives a large set of
  projective examples.
\end{remark}

%%%%%%%%%%%%%%%%%%%%%%%%%%%%%%%%%%%%%%%%
\section{$K$-classes of degenerations}
\label{sec:k-classes-degen}

We follow the notation of Section \ref{sec:one-param-families}.
Fix a projective dimension $n$, a Hilbert function $h$, 
an $\naturals$-graded locally free $A$-algebra $R$, and a
homomorphism ${A}[x_0,\ldots,x_n]\to R$ making $R$ a finite
${A}[x_0,\ldots,x_n]$-module. Let $R'$ be the integral closure
of $R$ in $R\tensor_A K$. 

Since $R$ and $R'$ are finite modules over $A[x_0,\ldots,x_n]$,
we see that $R/tR$ and $R'/tR'$ are finite modules over $k[x_0,\ldots,x_n]$,
and define elements of algebraic $K$-homology
of $k[x_0,\ldots,x_n]$. Since $R$ and $R'$ are both locally free and agree 
after inverting $t$, these two elements $[R/tR],[R'/tR']$ 
of $K$-homology coincide.  We now give a direct proof of this
$K$-equivalence, allowing us to strengthen the statement.

\begin{proposition}\label{prop:k-equivalence}
  Let $t$ be a non-zero-divisor in $R$, and let $R'$ stand between $R$ and
  its integral closure in $R[t^{-1}]$. Assume that $R'$ is finite over $R$
  (e.g. if $R$ is finitely generated over a DVR $A$).
  
  Then $\exists\ N>0$ such that $R/tR$ and $R'/tR'$ are
  $K$-equivalent modules over the ring $R/(t^N)$.
\end{proposition}

\begin{proof}
  Consider the short exact sequences of $R$-modules
  $$ 0 \to R/tR \to R'/tR \to R'/R \to 0 $$
  $$ 0 \to tR'/tR \to R'/tR \to R'/tR' \to 0 $$
  Since $t$ is not a zero divisor, the natural map $R'/R \to tR'/tR$
  is an isomorphism. So we get the $K$-equation
  $$ [R/tR] = [R'/tR] - [R'/R] = [R'/tR] - [tR'/tR] = [R'/tR']. $$

  By the assumptions on $R'$, there exists an $N$ such that
  $t^{N-1} R' \subset R$. Therefore $t^N$ annihilates all of these modules,
  so they are modules over the ring $R/(t^N)$, and the derivation of
  this $K$-equation holds there.
\end{proof}

The fact that $R/tR$ and $R'/tR'$ have the same Hilbert polynomial says only
that they define $K$-equivalent sheaves on $\bP^n$. Both sheaves are 
supported on thickenings of the same subvariety. The above Proposition
says that they are already $K$-equivalent on some larger thickening
of this same variety.

Note that there is a \emph{ring} homomorphism $R/(t) \to R'/(t)$.
Hence, given a family over $\Spec K$ of subvarieties of projective
space, there is a map from the limit branchvariety to the limit
subscheme inducing this $K$-equivalence.

\begin{example}
  Recall the colliding skew lines from Example \ref{ex:skewlines}.
  In this case, $\Proj R'/(t)$ is the two disjoint lines and
  $\Proj R/(t)$ has the lines crossing with an extra point embedding
  at the cross. If $\pi : \Proj R'/(t) \to \Proj R/(t)$ denotes the
  obvious collapse, then $\pi_*$ of the structure sheaf on the
  two lines is $K$-equivalent, but not isomorphic, to the structure
  sheaf of $\Proj R/(t)$.
\end{example}

The map $R/(t) \to R'/(t)$ was studied in \cite{Knutson05} 
(where one can find many more examples) in the
balanced normal cone context of Section \ref{sec:balanc-norm-cones}.
A principal result of that paper was that the corresponding map
from the balanced normal cone to the ordinary normal cone takes the
fundamental Chow class to the fundamental Chow class, a consequence
of this lemma.

%%%%%%%%%%%%%%%%%%%%%%%%%%%%%%%%%%%%%%%%
\section{The forest of a branchvariety}
\label{sec:rooted-forests}

Hartshorne proved \cite{Hartshorne2} that Hilbert schemes are connected. 
This is in some sense a negative result; it says that the only locally
constant invariants are the embedding dimension $n$ and the Hilbert
polynomial $h$. We already proved in Lemma \ref{lem:degree-sequence} that
for branchvarieties, the degree sequence $b$ is an additional such
invariant. In this section we develop a still finer invariant,
assuming the characteristic is $0$ or large enough.

\begin{lemma}\label{lem:connectedcomponents}
  Let $f:X\to \bP^n_S$ be a family of branchvarieties. Then
  the number of connected components (not irreducible components)
  of $X$ is locally constant.
\end{lemma}

\begin{proof}
  It is sufficient to assume that $S$ is a germ of a one-dimensional
  scheme, for example the $\Spec$ of a DVR. Let $X_{\eta}$ be the
  generic fiber.  Make a ramified base change $A'\supset A$ so that 
  the connected components of $X_{\eta}\otimes_{k_\eta} \overline k_{\eta}$ 
  are already defined over $A'$; call them $X_i$. 
  Then $\coprod X_i\to X$ is finite and they have the same generic
  fiber, so by the separatedness of $\B$
  (Corollary~\ref{cor:separated}) these two spaces coincide.
\end{proof}

As Example \ref{ex:skewlines} of the two colliding lines shows, this
is not an invariant for connected families of subschemes (it
only behaves semicontinuously).

\begin{assumption}\label{assume:char}
  Now let us fix a Hilbert polynomial $h(d)$ and a sequence of
  nonnegative integers $b=(b_0, \dotsc, b_{\deg h})$. For the rest of
  this section we will assume that our base scheme $C$ is defined over
  $\mb Z[1/(\max b_i)!]$. In other words, all prime numbers $p\le \max b_i$ 
  are invertible, and in particular any field over $C$ either has
  characteristic zero or $\chr k=p> \max b_i$. (For example, one can
  take $C$ itself to be $\Spec$ of a field $k$ whose characteristic is $0$ or
  some $p> \max b_i$.)
\end{assumption}

With this assumption on characteristic, every
generic plane section $X_i$ is reduced, hence is itself a branchvariety.

In this case, families of branchvarieties have some additional locally
constant invariants. For each connected component $X(j)$ of $X$, the
number of connected components of its general hyperplane section
$X(j)_1$ is locally constant, and we can continue by induction. 

To organize this induction, recall first the definition of a
\defn{rooted forest}, which is a graph with no cycles and a choice of
``root'' vertex in each connected component.  The vertices of a rooted
forest naturally form a ranked poset, where the rank of a vertex is
the length of the unique path to a root, and $v\geq w$ if $v$'s path to 
a root goes through $w$. Each component of a rooted forest, minus its root,
is itself naturally a rooted forest whose roots are the neighbors of the
old root. This allows rooted forests to be described inductively.

\newcommand\rk{{\rm rk}}

\begin{definition}
  Define the \defn{(labeled rooted) forest $\Forest(X)$ of a branchvariety} $X$
  %over an algebraically closed field $k$ 
  inductively as follows:
  \begin{enumerate}
  \item $\Forest(X) = \coprod_{X(j)} \Forest(X(j))$,
    where $\{X(j)\}$ are the connected components of $X$. (This includes
    the case $X=\emptyset$.)
  \item If $X$ is connected, then $\Forest(X)$ has one root (so, it is
    also connected).  The \defn{label} on the root is the integer
    $\chi(\mc O_X)$, the constant term in the Hilbert polynomial.
    Removing the root leaves $\Forest(X_1)$, where $X_1$ is a general
    hyperplane section of $X$.  (In order to have general enough
    hyperplanes, we work with the geometric fiber $X\otimes_k \bar k$.)
  \end{enumerate}
  So each vertex $v$ has two numbers associated: its rank $\rk(v)$ and
  its label $\chi(v)$. Also, if we have picked specific plane sections
  $X \supset X_1 \supset X_2 \supset \ldots$ of $X$, then we can 
  speak of the branchvariety associated to $v$, meaning the corresponding
  connected component of $X_{\rk(v)}$. 
  (By Assumption \ref{assume:char}, this is again a branchvariety.)
  In this way we can reduce arguments
  about a general vertex $v$ to the case that $v$ is a root, and the only one.
\end{definition}

\begin{proposition}\label{prop:forest}
%  With Assumption~\ref{assume:char},
%%% Already globally assumed
  $\Forest(X)$ is locally constant in families of branchvarieties.
  The Hilbert polynomial of $X$ and its degree sequence $(b_i)$ can be computed
  from $\Forest(X)$ by the formulae
  $$ h_X(d) = \sum_{v\in\Forest(X)} \chi(v) {d+\rk(v)-1 \choose  \rk(v)},
  \qquad
   b_i = \hbox{ the number of leaves of rank $i$}  $$
   where a \defn{leaf} is a maximal element of the poset $\Forest(X)$.
\end{proposition}

\begin{proof}
  The first part follows from the argument in Lemma 
  \ref{lem:connectedcomponents} applied to generic plane sections
  (which are again branchvarieties, by Assumption \ref{assume:char})
  and the invariance of the Hilbert polynomial in flat families.
  
  For the Hilbert polynomial formula, it is enough to check that the
  two sides agree at $d=0$ (which is obvious) and also after applying the 
  differencing operator $\Delta$ defined by $(\Delta g)(d) = g(d)-g(d-1)$.  
  The left hand side is $(\Delta h_X)(d) = h_{X_1}(d)$. The right hand is
  \begin{eqnarray*}
  \sum_{v\in\Forest(X)} \chi(v)  \Delta {d+\rk(v)-1 \choose  \rk(v)}
  &=& \sum_{v\in\Forest(X)} \chi(v) {d+(\rk(v)-1)-1 \choose  \rk(v)-1} \\
  &=& \sum_{v\in\Forest(X_1)} \chi(v) {d+\rk(v)-1 \choose  \rk(v)}
  \end{eqnarray*}
  which by induction is also $h_{X_1}(d)$.
  
  A leaf corresponds to an isolated point of a plane section.
  The degree of an irreducible component can be
  computed as the number of points in a generic plane section of
  complementary dimension, which implies the formula given for $b_i$.
\end{proof}

If we write the Hilbert polynomial and degree sequence associated to
a forest $F$ as $h(F)$ and $b(F)$, we can write
$$ \B^b_{h,Y} = \coprod_{F:\ b(F)=b, h(F)=h} \B_{F,Y} $$
where $\B_{F,Y}$ denotes the evident substack of $\B^b_{h,Y}$.
Note that this right-hand side is a finite union, since the left-hand
side has only finitely many connected components by Theorem \ref{thm:main}; 
only finitely many $F$ with $b(F)=b, h(F)=h$ give nonempty $\B_{F,Y}$.

\begin{example}
  If $X$ is an irreducible branchvariety of dimension $n$ and degree
  $d$, then $X$'s forest looks like a lone \defn{palm tree}, 
  with one vertex of each rank $<n$ and $d$ vertices of rank $n$.  
  More generally, $X$'s forest is a palm tree iff $X$ is
  equidimensional and connected in codimension $1$.
\end{example}

\begin{example}
  Consider two branchvarieties (indeed, reduced subschemes) of $\bP^3$:
  \begin{enumerate}
  \item $X=X(1)\sqcup X(2)$, where $X(1),X(2)$ are each a union 
    of two copunctal lines;
  \item $X'=X'(1)\sqcup X'(2)$, where $X'(1)$ is a line, and 
    $X'(2)$ is a union of three lines that pass though a common point but do
    not lie in a common plane.
  \end{enumerate}
  Then $X$ and $X'$ have the same Hilbert polynomial $p(d)=4d+2$ and
  the same degree sequence $b_0=0$, $b_1=4$. However, $\Forest(X)$
  consists of two trees with two leaves each, and $\Forest(X')$
  consists of a tree with three leaves and a tree with one leaf. $X$ and
  $X'$ therefore belong to different connected components of $\B$.
\end{example}

\begin{example}
  Consider two branchvarieties of $\bP^1$:
  \begin{enumerate}
  \item $X=X(1)\sqcup X(2)$, where each of $X(1),X(2)$ is a union of two 
    tangent $\PP^1$s;
  \item $X'=X'(1)\sqcup X'(2)$, where $X'(1)$ is a union of two
    crossing $\PP^1$s and $X'(2)$ is union of two $\PP^1$s meeting
    in a triple point.
  \end{enumerate}
  (Each irreducible component is degree $1$ as a branchvariety.)

  Then $X$ and $X'$ have the same Hilbert polynomial $p(d)=4d$,
  the same degree sequence $b_0=0$, $b_1=4$, and the same
  {\em unlabeled} rooted forest. However, the labels at the two roots
  of $F(X)$ are $0,0$ and at the two roots of $F(X')$ are $+1,-1$.  
  $X$ and $X'$ therefore belong to different connected components of $\B$.
\end{example}

\begin{example}\label{ex:stanley-reisner}
  If $F$ is a forest with all labels $1$, then $\Proj$ of the
  Stanley-Reisner ring of the order complex of the poset $F$ is a
  reduced subscheme of projective space whose forest is $F$.  So every
  (finite) labeled rooted forest is, up to relabeling, 
  the forest of some branchvariety.
\end{example}

However, not every labeling can occur; here is one of the simplest 
required conditions.
We give another necessary condition in Proposition \ref{prop:top-fork}.

\begin{proposition}\label{prop:only-one-leaf}
  Let $F$ be the forest of a branchvariety $X$, and $v$ a vertex of
  $F$ with only one leaf above it (not necessarily directly above).
  Then the label on $v$ is $1$.
\end{proposition}

\begin{proof}
  We can reduce to the case that $X$ is connected and $v$ is the root.
  Then by Proposition \ref{prop:forest}, $X$ is degree $1$, so by
  Zariski's Main Theorem the map $f:X\to \PP^n$ is the inclusion
  of a linear subspace. Hence $\chi(v) = \chi({\mc O}_{\PP^{\dim X}}) = 1$.
\end{proof}

\begin{question}\label{q:labelings}
  Let $f$ be a rooted unlabeled forest. Which labelings $F$ of $f$ 
  actually arise as forests of branchvarieties?
\end{question}

The corresponding question for subschemes (namely, which polynomials can arise 
as Hilbert polynomials) was solved by Macaulay (see \cite{Stanley,Green}).

\begin{question}\label{q:connected}
  For those $F$ that {\em are} forests of branchvarieties,
  is the open and closed subset $\B_{F,\PP^n}$ of $\B^{b(F)}_{h(F),\PP^n}$
  corresponding to a fixed forest $F$ connected?
\end{question}

As we will see in Corollary \ref{cor:connectivity} in the next
section, this connectivity depends only on $F$ and not the
dimension of $\PP^n$ (as long as $n$ is greater than or equal to the
maximum rank of $F$, for otherwise there are no branchvarieties).

\begin{example}
  Abandon Assumption~\ref{assume:char} on characteristics for this example.
  Assume $\chr k=p>0$ and 
  consider a family of branchvarieties $X=\bP^1$ of $Y=\bP^1$ parametrized by
  $S=\bA^1_t$, with $f$ given by the formula
  $$ (x_0,x_1) \mapsto (y_0=x_0^p,\ y_1 = x_1^p -t\cdot x_1 x_0^{p-1})$$
  For $t\ne0$ this is an unramified map, in fact Galois with Galois
  group $\mb Z/(p)$, and $\Forest(X)$ is a tree with $p$ leaves. 
  For $t=0$ we get a geometric Frobenius, a purely inseparable map,
  and $X_1$ is not a branchvariety.

  Hence, in this case the $p$ leaves ``glue together'' into a ``thick''
  branch. 
%It is possible that through such procedures, components for different forests
%may glue  together if  $0< \chr k \le \max b_i$. 
% I could not think of an example, and I don't want to say a wrong thing. So
% it is safest to skip this comment.
\end{example}

%%%%%%%%%%%%%%%%%%%%%%%%%%%%%%%%%%%%%%%%
\section{$U$-invariant branchvarieties}
\label{sec:u-invts}

In this section, we operate under the same Assumption~\ref{assume:char} 
on characteristics.

The space $\B_{F,\PP^n}$ of branchvarieties of $\PP^n$ with forest
$F$ carries an action of $\PGL_{n+1}$.
Let $U_{n+1}$, or just $U$, denote the group of upper triangular
matrices with $1$s on the diagonal. This acts on $\PP^n$ with $n+1$ orbits; 
two points are in the same orbit if they have the same last nonvanishing 
homogeneous coordinate. We will call the closures of these orbits the 
\defn{standard} $\PP^d$s in $\PP^n$. One motivation for studying 
actions of $U$ is the following:

\begin{lemma}\label{lem:u-connectivity} \cite{Horrocks}
  Let $U$ act on a complete space $X$, 
  and let $X^U$ be its fixed point set. 
  Then $X^U$ is connected iff $X$ is connected, and 
  the inclusion $X^U\to X$
  induces an isomorphism on
  (algebraic) fundamental groups.
\end{lemma}

(It is also easy to show that 
$X^U$ is rationally connected iff $X$ is,
which we will neither prove nor use.)

\junk{
\begin{proof}
  The property we use of $U$ is that it has a filtration by normal subgroups
  with all subquotients ${\mb G}_a$, the additive group. Using these
  one-parameter groups, it is easy to connect any point $x\in X$ to a
  $U$-fixed point by a chain of rational curves. 

  Conversely, let $x,y$ be two points in $X^U$, connected by a chain
  of (rational) curves. Again using the one-parameter groups in $U$,
  we can degenerate the chain to be \emph{setwise} $U$-invariant.
  Components can break under such a degeneration, but can't become nonrational.

  Since $U$ is connected, it acts on each irreducible component of this chain.
  If $U$ acts nontrivially on an irreducible curve, it can only fix
  one point, so we can excise such components without breaking connectivity.
  Hence $x,y$ are connected in $X^U$ by a chain of (rational) curves.
  Thus if $X$ is (rationally) connected, so is $X^U$.
\end{proof}
}

%\begin{remark}
%  We note that since $U$ is connected and the automorphism groups of are
%  finite, a branchvariety $f:X\to Y$ is $U$-invariant if and only if the
%  corresponding point in the coarse moduli space is. Hence, we can work
%  equivalently 
%  with the $U$-invariant points either in the stack or in the
%  coarse moduli space.
%\end{remark}

\begin{corollary}\label{cor:connectivity}
  Let $F$ be a rooted forest with maximum rank $d$. 
%  Let $\B_{F,\PP^n}$ denote the stack of branchvarieties
%  of $\PP^n$ with forest $F$.
%%% That was already defined in the previous section
  Then the number of components of $\B_{F,\PP^n}$ is constant for
  $n\geq d$ (and $0$ for $n<d$).

  More precisely, let $n_1\geq n_2\geq d$. Then the natural inclusion
  $\B_{F,\PP^{n_2}} \into \B_{F,\PP^{n_1}}$ induces isomorphisms
  on $\pi_0$ and $\pi_1$.
\end{corollary}

\begin{proof}
  The dimension of a branchvariety is the maximum rank of the vertices
  in its forest. Hence the $n<d$ statement is trivial; in this case
  there are no finite maps from a $d$-dimensional scheme to $\PP^n$.

  For a branchvariety $f:X\to\PP^n$ to be $U$-invariant (as a point of
  $\B_{F,\PP^n}$), its image must be $U$-invariant. The only $d$-dimensional
  $U$-invariant reduced subscheme of $\PP^n$ is the $\PP^d$ with
  vanishing last $n-d$ coordinates. Hence
  $$ (\B_{F,\PP^n})^{U_{n+1}} \cong (\B_{F,\PP^d})^{U_{d+1}}. $$
  Now apply the lemma.
\end{proof}

The case $n=d$ has the following classical description: it
parametrizes reduced schemes equipped with Noether normalizations.

We now attempt to better describe the elements of the $U$-fixed point set. 

\begin{lemma}\label{lem:irr-u-invts}
  Let $f:X^d\to\PP^n$ be an \emph{irreducible} branchvariety, defining
  a $U$-invariant point in the moduli stack. (In this case all
  Assumption \ref{assume:char} says is that $\chr k=0$ or $\chr k>\deg X$.) 
  Then $f$ is just the inclusion of the standard $\PP^d$ into $\PP^n$.
\end{lemma}

\begin{proof}
  In the proof of Corollary \ref{cor:connectivity}, we already determined 
  the image. So we may as well assume $d=n$ for ease of description.
  
  Let $\AA^n$ denote the open $U$-orbit, with last coordinate
  nonvanishing.  Then $X^\circ := f^{-1}(\AA^n)$ is open and nonempty
  in $X$, hence irreducible.  The map $f^\circ:X^\circ \to \AA^n$ is a
  $U$-invariant branchvariety of $\AA^n$. By the $U$-invariance,
  $f^\circ$ is ramified either everywhere or nowhere. 

% Terminology is fine here.
  
  If it is ramified everywhere (which can only happen in characteristic $p$), 
  then the degree of this cover is at least $\chr k$. But this violates
  our assumption on the characteristic.

%  With our assumption on the characteristic, $\mb A^n$ does not have
%  irreducible finite unramified covers of degree $\deg X$.
  Hence $f^\circ$ is a trivial cover, and for $X^\circ$ to be irreducible
  $f^\circ$ must be a degree $1$ cover,
  again by the assumption on the characteristic.
  So $f:X\to \PP^d$ is degree $1$.
  By Zariski's Main Theorem, $f$ is an isomorphism.
\end{proof}

\begin{example}
  We note that the conclusion of this Lemma is no longer true in small
  characteristics. Indeed, let $\chr k=p$ and let $f^{\circ}:\mb A^1\to \mb
  A^1$ be a morphism defined by $y=x^p-x$; homogenize to obtain a
  branchvariety $f:X=\bP^1\to Y=\bP^1$. Then $f^{\circ}$ is a nontrivial
  \'etale cover; in fact it is Galois with the Galois group $\mb Z/(p)$. One
  easily checks that the branchvariety $X$ is $U$-invariant.
\end{example}

\begin{lemma}\label{lem:u-eqvt}
  Let $f:X\to\PP^n$ be a branchvariety defining a $U$-invariant
  point in the moduli stack.
% and $\chr k=0$ or $\chr k>b_i$ for each  $b_i = \deg X^{\dim i}$. 
%%% Already assumed at the beginning of this section.
  Then there is an action of $U$ on $X$ such that $f$ is equivariant
  (and this action is unique).
\end{lemma}

\begin{proof}
  Let $J \subset \Aut(X,L) \times U$ be the closed subscheme of pairs $\{
  (a,u) : f\circ a = u\circ f \}$. (We note that $\Aut(X,L)$ is a closed
  subgroup of $\PGL_{N+1}$ and is an algebraic group.)  Its projection $p_2$ to
  the second factor is onto, by the assumption of $U$-invariance.
  
  The scheme $J$ carries an action of $\Aut(f)$ which on points is defined by
  $g\cdot (a,u) := (g\circ a,u)$, and it is easily seen to be free.  We note
  that $J$ is just the automorphism group scheme for a family of
  branchvarieties over $U$, as defined in Section
  \ref{sec:constr-moduli-space}. Hence $J\to U$ is finite. Since the action is
  free, $J\to U$ is \'etale.  By the assumption on the characteristic, $U$
  does not have irreducible finite \'etale covers of small degrees. So $J$
  breaks into a disjoint union of sections. Writing the section through
  $(1,1)$ as $u\mapsto (\alpha(u),u)$, the map $\alpha: U\to \Aut(X,L)$ is
  easily seen to give an action, and we are done.
\end{proof}

\begin{theorem}\label{thm:u-invts}
  Let $f:X\to\PP^n$ be a $U$-invariant branchvariety in $\B^F_{h,\PP^n}$.
%  where $\chr k$ is either $0$ or more than any $b_i(F)$. 
%%% This is a running assumption.
  Then $X$ is a union of projective spaces, where the
  number of irreducible components of dimension $i$ is
  the number of leaves of $F$ of rank $i$. Any (nontrivial)
  intersection of components is irreducible, and identified by $f$
  with a $U$-invariant thickening of one of the standard $\PP^d$s in
  $\PP^n$. 
\end{theorem}

\begin{proof}
  Let $X(1),\ldots,X(m)$ be a nonempty set of components of $X$, 
  and $H = \cap_j X(j)$, considered as a subscheme of $X(1)$. Then under 
  the identification $f:X(1) \wt\longrightarrow \PP^{\dim X(1)}$ guaranteed
  by Lemma \ref{lem:irr-u-invts}, $H\red$ maps to a reduced $U$-invariant
  subscheme of the standard $\PP^{\dim X(1)}$. The only possibility is
  the standard $\PP^{\dim H}$.
  
  Since the components are all degree $1$, the number of components
  of dimension $i$ is the degree of $X^{\dim i}$, which we already knew
  to be the number of leaves of $F$ of rank $i$.
\end{proof}

\begin{remark}
  It is not difficult to see that the seminormalization of any
  $U$-invariant branchvariety is exactly the Stanley-Reisner scheme from
  Example \ref{ex:stanley-reisner} for the same rooted forest (but all
  labels changed to $1$).
\end{remark}

\begin{proposition}\label{prop:top-fork}
  Let $F$ be the forest of a branchvariety $X$, and $v$ a maximal fork
  of $F$, i.e. each $w>v$ has only one leaf above $w$.
  Then the label $\chi(v)$ on $v$ is at most $1$. 
\end{proposition}

\begin{proof}
  As already explained, we can reduce to the case that $X$ is connected
  and $v$ is the root.
  By Lemma \ref{lem:u-connectivity} and Proposition \ref{prop:forest},
  we can assume $X$ is $U$-invariant. Then the condition on $v$, 
  plus the proof of Proposition \ref{prop:only-one-leaf}, says that the
  seminormalization $\wt X$ is an isomorphism away from the standard $\PP^0$.
  If $Y$ is the scheme-theoretic fiber lying over $\PP^0$ (a fat point), then 
  $\chi({\mc O}_X) - \chi({\mc O}_{\PP^0}) 
  = \chi({\mc O}_{\wt X}) - \chi({\mc O}_Y) = 1 - len(Y)$, where
  $len(Y)\geq 1$ is the length of the fat point $Y$.  
  Hence $0 \geq \chi({\mc O}_X) - \chi({\mc O}_{\PP^0}) = \chi(v)-1$.
\end{proof}
  
(We won't need it, but even if the $X$ in the above proof is not
assumed $U$-invariant, one can give a good description of it: $X$ is the 
connected union of a set $\{P_i\}$ of projective spaces, each including into
$\PP^n$ as a linear subspace, and a set of curves $\{C_j\}$, where
all intersections $P_i \cap P_j$ are $0$-dimensional.)

\begin{example}
  Not all labels on forests of branchvarieties are at most $1$.
  If $X$ is a quartic hypersurface in $\PP^3$, e.g. a K3 surface,
  then the label on the root of its forest is $2$. 
\end{example}

%%%%%%%%%%%%%%%%%%%%
\subsection{Spaces of cubic curves}
\label{sec:complete-cubics}

Consider the space of branchvarieties of $\PP^2$ including the plane cubics,
so $h(n) = 3n$, $b_0 = 0$, $b_1 = 3$. According to the
description in Theorem \ref{thm:u-invts}, a $U$-invariant branchvariety 
is a union of three copunctal lines. These already appeared in
Example \ref{ex:3lines}, where the angle of intersection gave a
$\PP^1/S_3$ worth of such branchvarieties.

On the other hand, there is a $3$-dimensional space of $U$-invariant 
stable curves of degree $3$ and genus $0$. Each has
an elliptic curve (which collapses entirely) meeting three $\PP^1$s
each in a point. So the space of branchcubics does not match the
corresponding moduli space of stable maps.

%%%%%%%%%%%%%%%%%%%%%%%%%%%%%%%%%%%%%%%%
\section{Relations to other moduli spaces}
\label{sec:conn-hilb-scheme}

In each of $\Hilb$, $\B$, and $\Chow$ there is an open set corresponding to 
reduced subschemes, and these three open sets are naturally identified.
In general, the only natural morphisms extending this identification
go from $\Hilb$ or $\B$ to $\Chow$, as the following examples show.

%%%%%%%%%%%%%%%%%%%%%%%%%%%%%%%%%%%
\subsection{Branch vs. Hilbert}
\label{sec:branch-vs-hilbert}

\begin{example}
  The moduli stack of $n$ points branched over a reduced scheme $Y$ is
  easily computed to be the ``symmetric product'' stack $[Y^n / S_n]$.
  If $Y = \PP^2$ and $n>1$ then the Hilbert scheme of $n$ points has a
  nontrivial blowdown to the \emph{coarse} moduli space $Y^n/S_n$
  (which is in fact the Chow variety).
  However, there is no natural \emph{stack} map $\Hilb \to \B$, and no 
  continuous map $\B \to \Hilb$.
\end{example}

\begin{example}
  The Hilbert scheme of plane conics is simply $\PP^5$. The branchvariety
  stack of plane conics is (coarsely) $\PP^5$ blown up along the Veronese 
  surface. So in this case there is no continuous map $\Hilb \to \B$.
\end{example}

On the other hand, as was pointed out to us by J\'anos Koll\'ar, the classical
Hilbert scheme can be properly embedded into $\B$,
albeit for different parameters.

\begin{theorem}\label{thm:Kollar}
  There exists a closed embedding 
  $\psi:\Hilb_{h,\bP^n} \to \B^{(0,\ldots,0,2)}_{2q-h,\bP^n}$,
  where $q(d) = \chi(\mc O_{\bP^n}(d))$. Moreover, if $\deg h\le n-3$ then
  $\Hilb_{h,\bP^n}$ is (coarsely) a connected component of 
  $\B^{(0,\ldots,0,2)}_{2q-h,\bP^n}$. As a stack, 
  this connected component of $\B^{(0,\ldots,0,2)}_{2q-h,\bP^n}$ 
  is $\Hilb_{h,\bP^n}$ modulo a trivial ${\mb Z}/(2)$ action.

  By Proposition \ref{prop:forest}, each branchvariety in this connected 
  component has the same forest $F$. If $\chr k\neq 2$, then
  the converse holds; each branchvariety
  with forest $F$ is of the sort just constructed. In particular,
  $\B_{F,\PP^n}$ is connected, giving a positive answer to
  Question \ref{q:connected} for this $F$.
\end{theorem}

\begin{proof}
  The morphism $\psi$ associates to each family of subschemes
  $Z\subset \bP^n_S$ the following family $X$ of branchvarieties of
  $\bP^n$: $X$ consists of two copies of $\bP^n_S$ glued along $Z$;
  the structure sheaf  of $X$ is 
  $\varprojlim(\mc O_{\bP^n_S}\oplus \mc O_{\bP^n_S}\to \cO_Z)$.
  Clearly, this gives a closed embedding. 
  
  By \cite[Cor. 12.7]{Kollar_FlatnessCriteria}, if $\codim Z\ge 3$ then any
  deformation of $X$ extends to a deformation of the
  two copies of $\bP^n$, i.e. it is of the same type. In this case, the whole
  connected component consists of branchvarieties glued from two copies of
  $\bP^n$, and $Z$ can be recovered from it.  The two copies of $\bP^n_S$ can
  be unambiguously called the ``first'' and ``second'' copy, without monodromy
  in $Z$, hence the triviality of the ${\mb Z}/(2)$ action.
  
  It is easy to describe the forest $F$ of any one of these branchvarieties.
  Taking a plane section with a plane of codimension $k$ gives two
  copies of $\bP^{n-k}_S$ glued along a plane section of $Z$. 
  If $k\leq\deg h$, the plane section of $Z$ is nonempty, so the
  branchvariety $X_k$ is connected. For $k>\deg h$, the plane section
  $X_k$ is a union of two copies of $\PP^{n-k}_S$. So $F$ looks like a
  tuning fork; it has one vertex $v$ for each rank $0,\ldots,\deg h$,
  labeled $1-(\Delta^{\rk v} h)(0)$, the fork vertex at rank $\deg h$, 
  and two vertices (each labeled $1$, as to be expected from
  Proposition \ref{prop:only-one-leaf}) for each rank $k+1,\ldots,n$.

  Now we want to show that the only connected component of $\B^{(0,\ldots,0,2)}_{2q-h,\PP^n}$
  with this forest $F$ is the one described above. By Lemma 
  \ref{lem:u-connectivity}, it is enough to check that any
  {\em $U$-invariant} branchvariety $X$ with this forest is in
  the connected component above. But by Theorem \ref{thm:u-invts}
  (for which we need $\chr k\neq 2$),
  such an $X$ is obviously a union of two copies of $\PP^n_S$ along a 
  $(\deg h)$-dimensional scheme.
\end{proof}

\begin{example}
  Consider the case when $Z$ is two points (or one point of
  multiplicity 2) in $\PP^1$. In this case the Hilbert scheme is 
  $(\PP^1)^2/S_2 \cong \PP^2$. Whereas $\B$ is $(\PP^1)^4/S_4$,
  as explained in Example \ref{ex:double-covers}.
\end{example}

\begin{remark}
  This example was the case of $\codim Z=1$. When $\codim Z=2$, one can say
  that $(\Hilb_{2q-h,\bP^n})\red$ is a connected component of $(\B)\red$,
  since a degeneration of a branchvariety connected in codimension $1$ is again
  connected in codimension $1$. However, the scheme structures in this case
  are not obvious.
\end{remark}

%%%%%%%%%%%%%%%%%%%%%%%%%%%%%%%%%%%%%%%%
\subsection{Branch vs. Chow}
\label{sec:branch-vs-chow}

The Chow variety of plane conics, like the Hilbert scheme, is $\PP^5$.
So there is no continuous morphism in general from Chow to $\B$.

If $f:X\to Y$ is a branchvariety of $Y$ then $(f_*[X^{\dim i}])$ is a
well-defined collection of cycles on $Y$, of dimension $i$ and degree $b_i$.
This gives a set-theoretic map from $\B$ to a product of Chow varieties.
Putting this into arbitrary families is somewhat delicate since, as we already
noted in~\ref{rem:Branch-vs-Chow}, Chow lacks the infinitesimal theory.
J. Koll\'ar \cite[I.3-4]{Kollar_RationalCurves} defines the Chow functor on the
category of reduced and seminormal schemes. This immediately gives a
functorial morphism
$$ (\Bhb)_{\red}^{\rm semi} \to \prod_i \Chow_{i,b_i,Y} $$
from the seminormalization of the reduced part of $\B$. 

We note that the information encoded in Branch is much richer, and the
cycle $f_*[X]$ is but a shadow of~$X$.

\subsection{Branch vs. stable maps}
\label{sec:branch-vs-stab}

Any proper map $X\to Y$ (for $X$ reduced) admits a Stein factorization
$X \to X' \to Y$ where $X'\to Y$ is a branchvariety. This suggests that
there might be natural transformations to $\B$ from other spaces of maps --
in particular, stable maps of curves.

However, $\chi(X,L^d) \neq \chi(X',L^d)$ in general, which would make
such a transformation discontinuous. One case in which they are equal
is $X\to Y$ a stable curve of genus zero, and indeed the Stein
factorization gives a natural transformation 
$\overline{M}_{0,0}(Y) \to \B_Y$.
The details will appear in \cite{Lin}.

%%%%%%%%%%%%%%%%%%%%%%%%%%%%%%%%%%%%%%%%
\section{Other versions}
\label{sec:other-versions}

%%%%%%%%%%%%%%%%%%%%%%%%%%%%%%%%%%%%%%%%
\subsection{Multigraded Branch}
\label{sec:mult-affine-branchv}

There are two multigraded analogues of the classical Hilbert scheme:
the toric Hilbert scheme \cite{PeevaStillman} and the multigraded
Hilbert scheme of \cite{HaimanSturmfels}. As was explained in
Section~\ref{sec:stable-toric-vari}, the toric $\B$ already exists: 
it is the moduli space of multiplicity-free stable toric varieties. 
To construct the multigraded $\B$ in full generality using the methods
of this paper would appear to require a definition of $b$-sheaf in the
equivariant setting, where invariant hyperplanes are not generic.

%%%%%%%%%%%%%%%%%%%%%%%%%%%%%%%%%%%%%%%%
\subsection{Complex-analytic Branch}

The complex-analytic analogues of the Hilbert scheme and Chow
varieties classifying complex-analytic subspaces, resp. cycles of a
complex-analytic space are well known; they are called the Douady
space and Barlet space respectively. Clearly, a complex-analytic
analogue of $\B$ can and should be constructed as well.

\bibliographystyle{amsalpha}
%\bibliography{bib-hilblike}

\providecommand{\bysame}{\leavevmode\hbox to3em{\hrulefill}\thinspace}
\providecommand{\MR}{\relax\ifhmode\unskip\space\fi MR }
\providecommand{\MRhref}[2]{%
  \href{http://www.ams.org/mathscinet-getitem?mr=#1}{#2}
}
\providecommand{\href}[2]{#2}

%\bibliography{primary}

\end{document}